\documentclass[11pt]{article}
\usepackage[active]{srcltx}
\usepackage{epsfig,amssymb}

\begin{document}

\title{
Local classification of singular hexagonal 3-webs\\ with holomorphic Chern connection and infinitesimal symmetries}

\author{{\Large Sergey I. Agafonov}\\
Department of Mathematics\\
Federal University of Paraiba,\\ Jo\~ao Pessoa, Brazil\\
e-mail: {\tt sergey.agafonov@mat.ufpb.br} }
\date{}
\maketitle
\unitlength=1mm

\newtheorem{theorem}{Theorem}
\newtheorem{proposition}{Proposition}
\newtheorem{lemma}{Lemma}
\newtheorem{corollary}{Corollary}
\newtheorem{definition}{Definition}
\newtheorem{example}{Example}

\pagestyle{plain}

\begin{abstract}
\noindent We provide a complete classification of  hexagonal singular 3-web germs in the complex plane, satisfying the following two conditions:

1)	the Chern connection remains holomorphic at the singular point,

2)	the web admits at least one infinitesimal symmetry at this point.\\
As a by-product, a classification of hexagonal weighted homogeneous 3-webs is obtained.
(secondary).
\end{abstract}

\section{Introduction}

A finite number of pairwise different foliations in the plane form a planar web. A point $q_0$ is called {\bf regular} if for each pair of foliations the tangent lines to the leaves at this point are transverse to each other. The corresponding local object is called a {\bf non-singular web germ} at $q_0$.
Consider the group of diffeomorphism (or biholomorphism)  germs and the corresponding equivalence relation. Clearly, each non-singular 2-web germ is equivalent to the 2-web germ formed by coordinate lines. The situation becomes more complicated for 3-webs.  Blaschke discovered that generically even a regular 3-web germ is not equivalent to the web germ of three families of parallel lines \cite{BB}. From the differential-geometric point of view a nontrivial 3-web has a non-vanishing curvature 2-form.

\begin{figure}[th]
\begin{center}

\epsfig{file=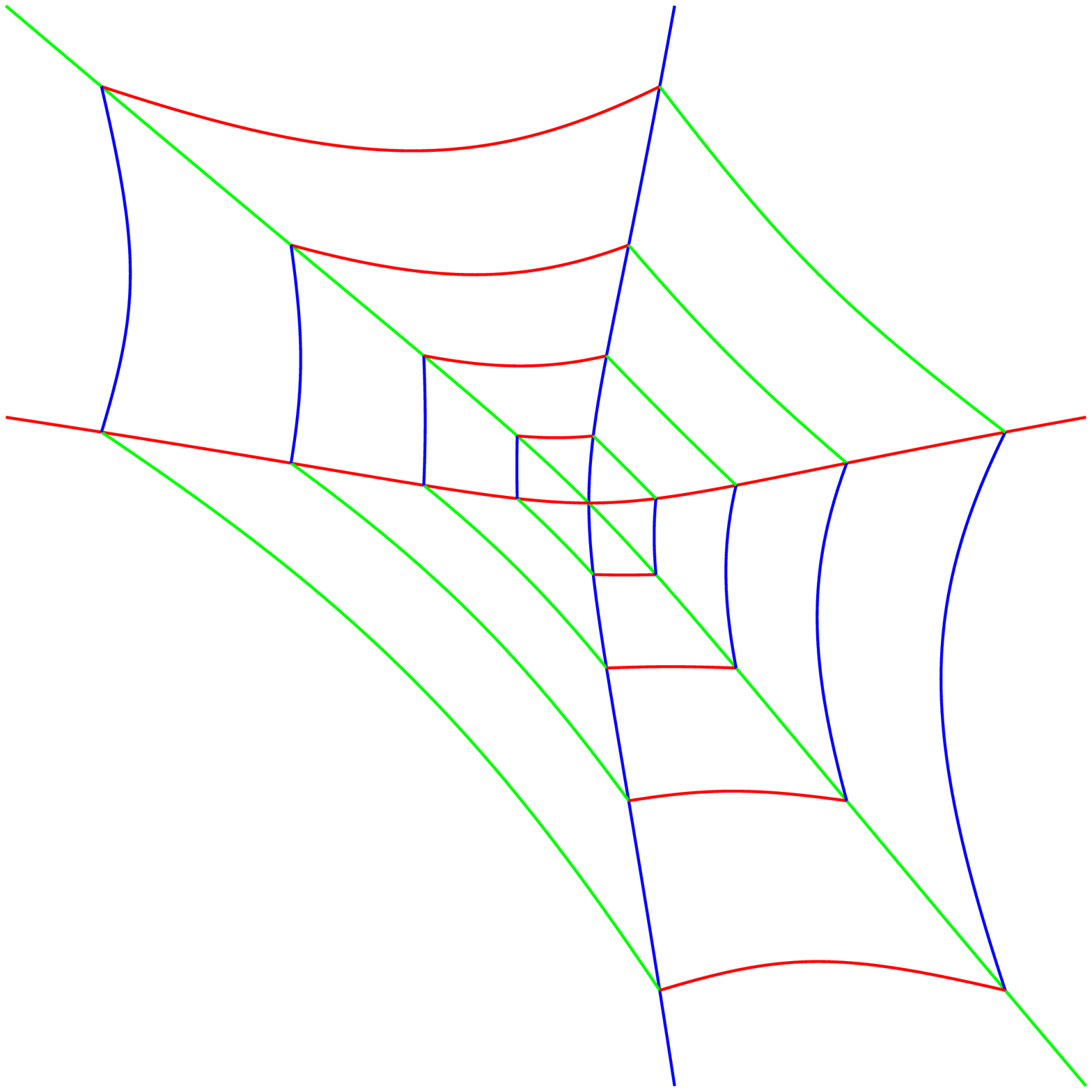,width=40mm}\put(-4,27){\large  $L_1$} \put(-21,38){\large  $L_2$} \put(-42,35){\large  $L_3$} \put(-19,18.5){ $\scriptstyle P_0$} \put(-10,20){ $\scriptstyle P_1$} \put(-19.5,20.5){\small $\bullet$} \put(-10.5,22){\small $\bullet$} \ \ \ \ \ \ \ \ \ \ \ \put(8,20){\vector(1,0){10}} \put(11,25){\huge  $\varphi$}  \ \ \ \ \ \ \ \ \ \ \ \ \ \ \ \ \ \  \ \ \ \epsfig{file=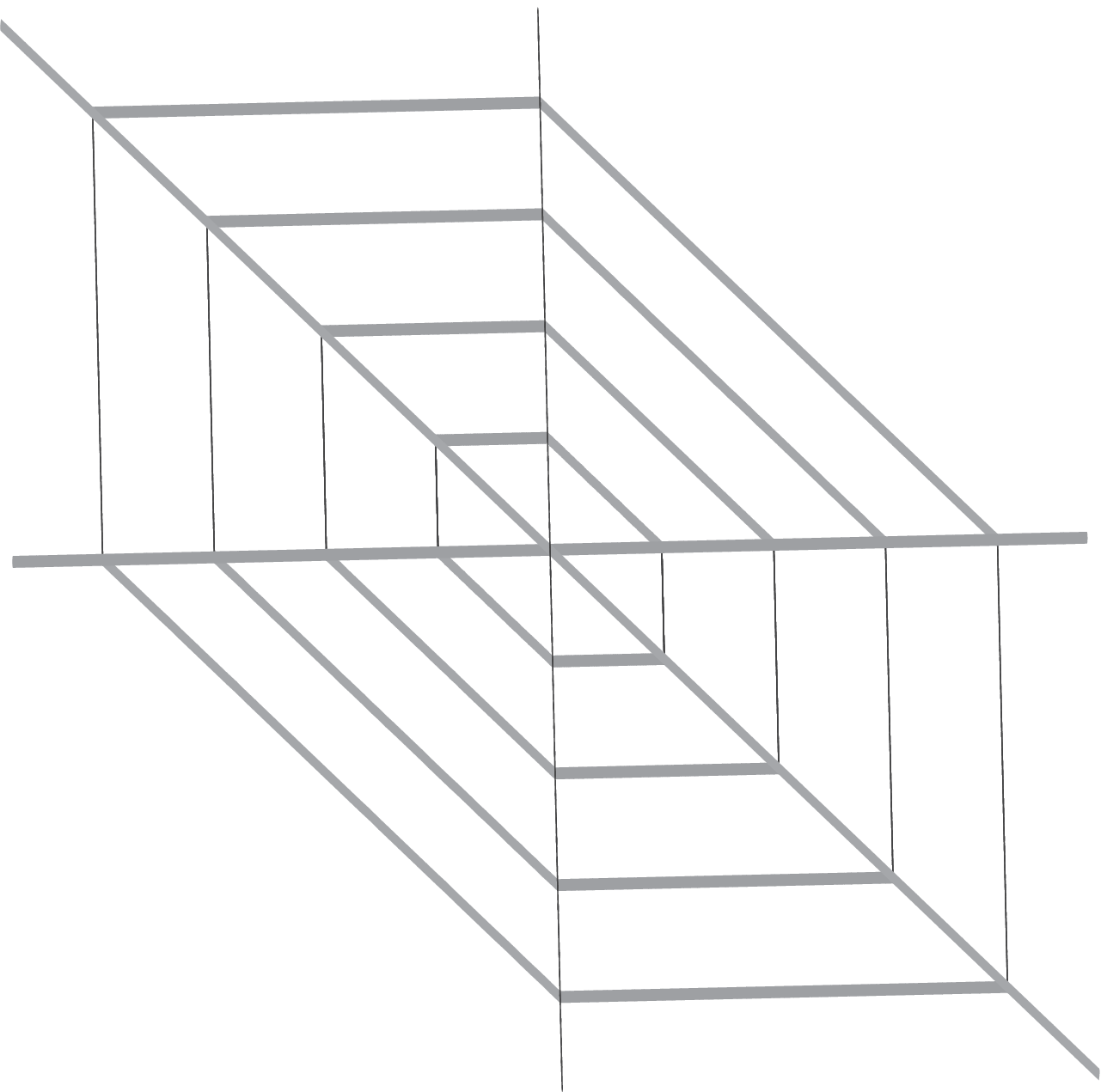,width=40mm}
  \caption{Brian\c{c}on's hexagons} \label{Pic_web}
\end{center}
\end{figure}

Moreover, the invariant in question is topological in nature. Choose an arbitrary regular point $p_0$ in the plane, draw the leaves $L_1,L_2,L_3$ of the web through this point, take a point $p_1$ on $L_1$ and go around $p_0$ along the web leaves starting from $p_1$ and swapping the foliation each time when meeting $L_1$, $L_2$ or $L_3$. For the web equivalent to 3 families of parallel lines, which has zero curvature, one comes back to $p_1$. Less trivial is that the inverse is also true and the following 3 conditions are equivalent:
 \begin{enumerate}
 \item for each choice of $p_0,p_1$, the constructed hexagon-like figure is closed,
 \item the web is equivalent to 3 families of parallel lines,
 \item the web curvature is zero.
 \end{enumerate}
 (see \cite{BB,Be} and picture \ref{Pic_web}).
The webs possessing any of the above properties are called {\bf hexagonal} or {\bf flat}.
A sufficiently general class of 3-web germs can be described by binary forms:
\begin{equation}\label{binary}
K_3(x,y)dy^3+K_2(x,y)dy^2dx+K_1(x,y)dydx^2+K_0(x,y)dx^3=0,
\end{equation}
where at least one of the germs $K_i$ does not vanish at $(0,0)$.  Dividing the above form by $dx^3$ and by a non-vanishing coefficient one gets an implicit ODE, cubic in $p=\frac{dy}{dx}$. Rotating the coordinate axes, if necessary, one can reduce this equation to a monic one:
\begin{equation}\label{general_cub}
p^3+a(x,y)p^2+b(x,y)p+c(x,y)=0.
\end{equation}
Its solutions form a hexagonal 3-web iff the coefficients of this cubic ODE satisfy a certain nonlinear partial differential equation (see Section \ref{Chern connection}).
In this paper we study the complex analytic case, i.e., $K_i$ are germs of holomorphic function at $(\mathbb C^2,0)$ and the equivalence relation is induced by the group of germs of biholomorphisms ${\rm Diff(\mathbb C^2,0)}$.

\begin{figure}[th]
\begin{center}
\epsfig{file=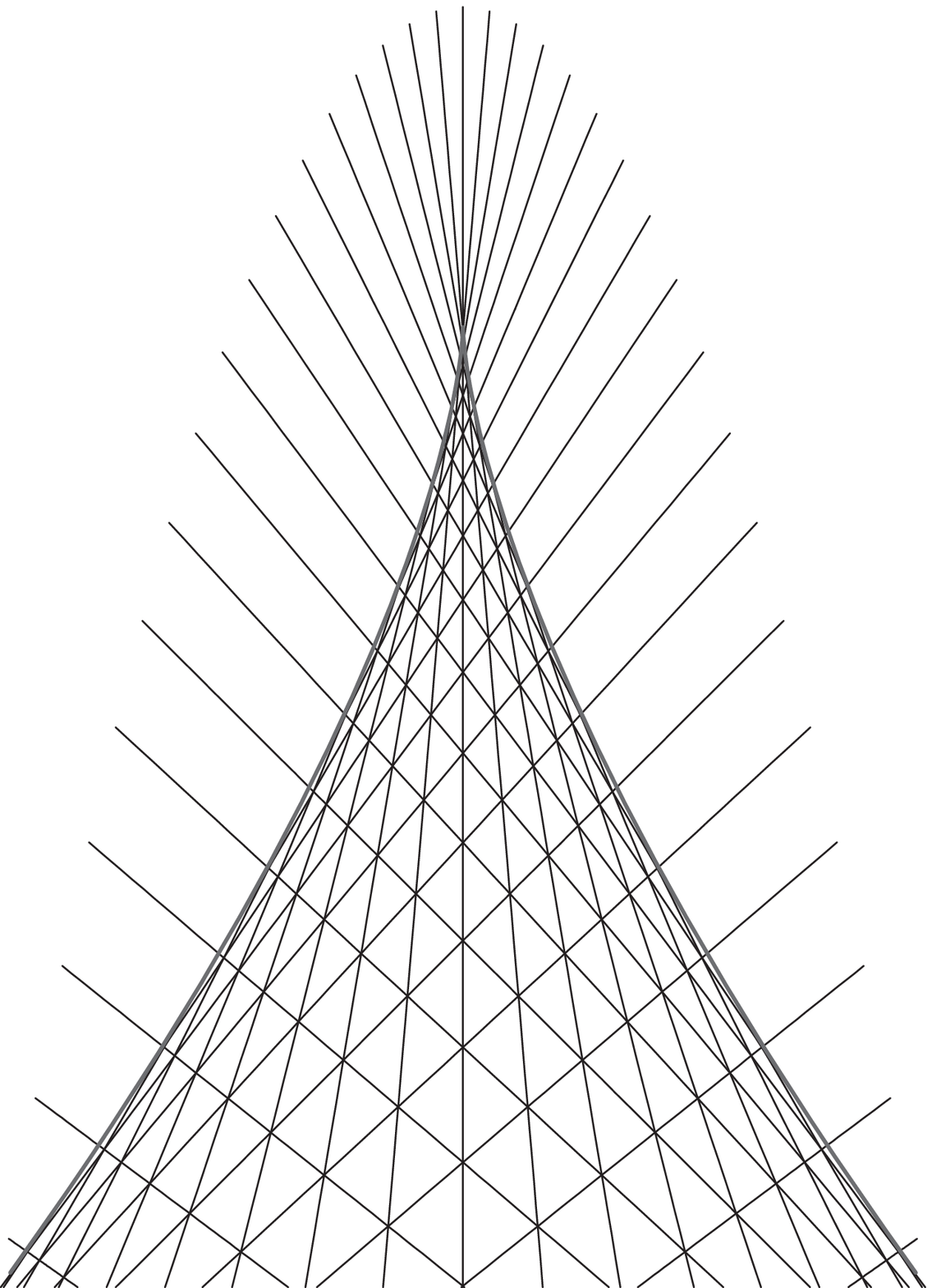,width=40mm} \ \ \ \ \ \ \ \ \ \ \ \ \ \ \ \
\ \ \ \ \ \ \ \epsfig{file=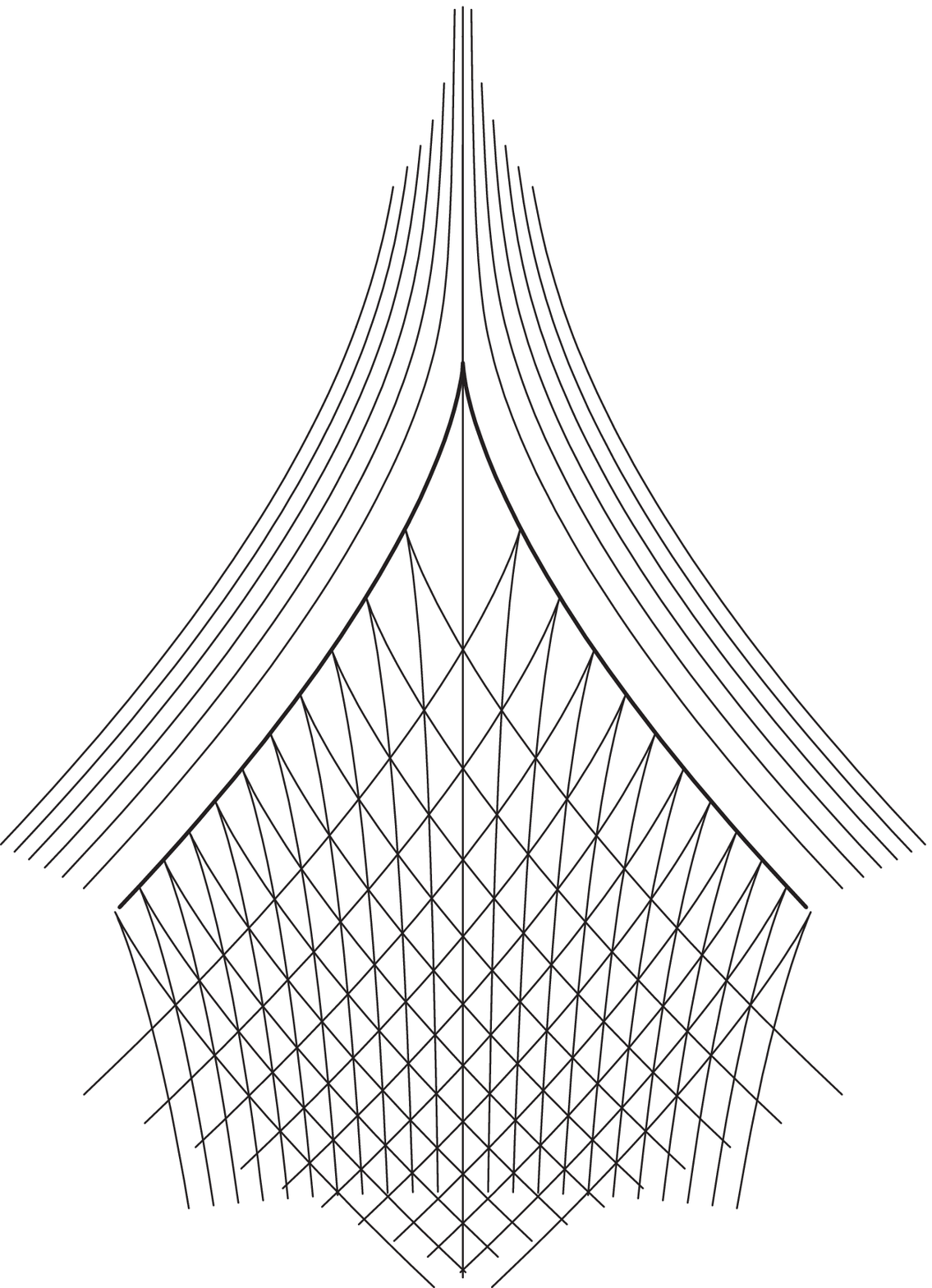,width=40mm} \caption{Solutions
of $p^3+px-y=0$ and $p^3+2xp+y=0$ with horizontal
y-axis.} \label{Pic_sol}
\end{center}
\end{figure}

\begin{example} {\rm The classical Graf and Sauer theorem
\cite{GSg} claims that a 3-web of straight lines is hexagonal iff
the web lines are tangents to an algebraic curve of class 3, i.e.,
the dual curve is cubic. This implies immediately that the
following cubic Clairaut equation has a hexagonal 3-web of
solutions:
$$
p^3+px-y=0.
$$
Its solutions are the lines $p=const$ enveloping a semicubic
parabola (Fig. \ref{Pic_sol}). }
\end{example}

\begin{example}\label{associativ}  {\rm  Consider an associativity equation
\begin{equation}\label{ass}
u_{xxx}=u^2_{xyy}-u_{xxy}u_{yyy},
\end{equation}
describing 3-dimensional Frobenius manifolds (see  \cite{Df}).
Each of its solutions $u(x,y)$ defines a characteristic web in the
plane. This web is hexagonal; it follows from the results obtained in \cite{MFa}. (See also \cite{Al}, where a broader class of PDEs with flat characteristic 3-webs was studied.)
Characteristics are integral curves of the vector field
$$
\partial _x-\lambda (x,y)\partial_y,
$$
where $\lambda$  satisfy the  characteristic equation
$$
\lambda ^3+u_{yyy}\lambda ^2-2u_{xyy}\lambda + u_{xxy}=0.
$$
For the solution $u=\frac{x^2y^2}{4}+\frac{x^5}{60}$,  the
characteristic equation becomes
$$
p^3+2xp+y=0
$$
after the substitution $x \to -x$, $y\to -y$, $\lambda\to -p$.  The corresponding 3-web is shown in Fig. \ref{Pic_sol}. }
\end{example}
Notice that the above web germs are not equivalent. We call a web germ at $q_0\in \mathbb C^2$ {\bf singular} if at least two web directions coincide at $q_0$. The examples show that singular hexagonal web germs are not necessarily equivalent.

Curvature 2-form of a 3-web is defined as the external derivative $d(\gamma)$ of the Chern connection 1-form $\gamma$ (see \cite{Be} and Section \ref{Chern connection}). Thus, for hexagonal 3-webs, this form is closed. But it is not exact in general: on the {\bf discriminant curve} of the web, where at least 2 foliations of the web are not transverse, the Chern connection form usually has a pole. For instance, for the first of the above examples we have
$$
\gamma=\frac{6 x^2 dx+27ydy}{4 x^3+27 y^2},
$$
whereas in the second example $\gamma$ is the zero 1-form, i.e. holomorphic. We are particularly interested in classification of singular webs whose Chern connection form remains holomorphic on the discriminant curve.

Observe that the above two 3-webs are invariant under the flow of the vector field $X=2x\partial _x+3y\partial _y.$
We say that web has an infinitesimal symmetry $$
X=\xi(x,y)\partial _x+\eta (x,y)\partial _y,
$$
if the local flow  of the vector field $X$ maps the web leaves to the web leaves. Infinitesimal symmetries form a Lie algebra with respect to the
Lie bracket. Cartan  proved  (see \cite{Cg}) that at a regular point
a 3-web either does not have infinitesimal symmetries
(generic case), or has one-dimensional symmetry algebra (then in
suitable coordinates it can be defined by the form
$dx\cdot dy \cdot (dy+u(x+y)dx)=0$ with the symmetry $\partial _y
-\partial _x$), or has a three-dimensional symmetry algebra (then
it is equivalent to the web defined by the form  $dx\cdot dy \cdot (dy+dx)=0$ with the
symmetry algebra generated by  $\{\partial _x,\
\partial _y, x\partial_x+y\partial_y\}$).
In the last case, when the symmetry algebra has the largest possible  dimension 3,
the 3-web  is  hexagonal. Note that not all symmetries survive at a singular point; in the above examples the dimension of
the symmetry algebra drops to 2 at a generic point of the discriminant curve and to 1 at the
cusp point. The condition to have at least one-dimensional
symmetry at a singular point is not trivial. The following equation has a flat 3-web of solutions but
does not admit non-trivial symmetries at $(0,0)$
 $$p^3-2x^2y(1+x^2)p+8x^3y^2=0.$$

\noindent The objective  of this paper is {\it to describe singular hexagonal 3-web germs such that:
\begin{enumerate}
\item the Chern connection form is holomorphic,
\item the web has at least one infinitesimal symmetry.
\end{enumerate}}

\noindent The principal motivation for this problem comes from the geometric theory of Frobenius manifolds. Namely characteristics on solutions of WDVV associativity equation (see
Example \ref{associativ}) form a  hexagonal 3-web, as was observed by Ferapontov. In fact, for this web, the Chern connection is holomorphic,  i.e locally exact.
Moreover, the associativity condition in suitable flat coordinates assumes two essentially different forms: either as equation (\ref{ass}) or as the following one (see  \cite{Df}):
\begin{equation}\label{ass2}
u_{xxx}u_{yyy}-u_{xxy}u_{xyy}=1.
\end{equation}
Now the characteristic 3-webs are defined by the following ODE:
$$
u_{yyy}p^3+u_{xyy}p^2-u_{xxy}p-u_{xxx}=0.
$$
We observe that the Chern connection is zero {\it on the characteristic 3-web of each solution $u(x,y)$} of equation (\ref{ass}) and is holomorphic {\it on the characteristic 3-web of each solution} of (\ref{ass2}). See \cite{Aw} and \cite{Af} for more detail and a geometrical interpretation of this fact. Moreover, the characteristic 3-web for a Frobenius 3-web germ can be constructed in a pure geometrical way sturting from a Frobenius 3-fold, as was shown in \cite{Aw}.

Further, if a solution of the associativity equation corresponds to some geometric Frobenius structure, as defined by Dubrovin in \cite{Df}, the characteristic 3-web has an infinitesimal symmetry, inherited from the so-called {\bf Euler} vector field (see \cite{Aw})
\begin{equation}\label{euler}
E=w_x x\frac{\partial}{\partial x}+w_y y\frac{\partial}{\partial
y}, \ \ \ \ \ \ w_x,w_y={\rm const}.
\end{equation}
Equations and webs symmetric with respect to such dilatation symmetry are called {\bf weighted homogeneous}. In what follows we call the operator of type (\ref{euler}) Euler vector field.

Hexagonal 3-webs, satisfying the above two conditions, have nice properties from the purely mathematical point of view. Recall that a 3-web is hexagonal iff its foliations have first integrals $u_i$ satisfying $u_1+u_2+u_3=0$. Then the finiteness of the Chern connection form implies that these integrals are integer algebraic over the ring of holomorphic function germs.
The corresponding algebraic equation for $u_i$ at a singular point $q_0$ can be read off the infinitesimal symmetry,  provided that this symmetry has an equilibrium point at $q_0$.
Further, if the monodromy group, permuting the web leaves on going around the discriminant curve, is "maximal possible", i.e. $S_3$ for the triple root of equation (\ref{binary}) or $Z_2$ for a double root, we prove also the existence of infinitesimal symmetries. Namely we prove that if the Chern connection form is exact and defined at some neighborhood of
a singular point $q_0$, then there is at least 2-dimensional symmetry algebra at $q_0$ for a double root and
at least 1-dimensional symmetry algebra for a triple root. (See section \ref{symsection}.)

The main result of the paper is a complete classification of the 3-web germs of the class introduced above. The list consists of 5 equations and 3 infinite series (see Theorems  \ref{homoexact} and \ref{shiftTH}). It is remarkable that the normal forms can be written in terms of polynomials, the function $\tan$ and Legendre`s functions $P^{\mu}_{\nu},Q^{\mu}_{\nu}$.   The key observation that allowed the obtained classification is that each  symmetry operator vanishing at the singular point is equivalent to dilataion symmetry (see Theorem \ref{Xeuler}).
Thus, as  a by-product, we obtained a complete classification of weighted homogeneous hexagonal 3-webs, i.e. having an infinitesimal symmetry of the form (\ref{euler}) (and possibly singular Chern connection).  We also provide invariants distinguishing the normal forms for these two classifications.

Symmetries, vanishing at singular points,  comes quite natural by singular webs.
For instance, it is immediate that if equation (\ref{general_cub}) admits a
non-trivial symmetry at a singular point of the discriminant curve, this point is necessary
a singular point of the vector field $X$.

In the literature mainly symmetries of
explicit ODEs were  studied. As the considerations were local this
is equivalent to the case of regular points, where the equation
can be resolved with respect to the derivative $p$. See
\cite{Lds},\cite{Lki} for classical treatment and \cite{Oa} for a modern
exposition.
If the ODE is explicit or quadratic  with respect to the derivative then
its symmetry algebra is infinite dimensional at a generic point. Indeed, an explicit equation can be brought to the form $dy=0$ and each operator of the form $\xi (x,y)\partial _x+\eta (y) \partial _y$ is a symmetry. For a quadratic ODE at a point with distinct roots we can choose the first integrals as local coordinates $x,y$, then ODE takes the form $dy\cdot d x=0$ and the symmetries are $\xi(x)\partial _x+\eta (y) \partial _y$. Thus the computing of infinitesimal symmetries of ODEs is related to integrating of  PDEs, whose solutions involve arbitrary functions. As this can not be done in general, the case of ODEs of the first order was rarely studied in the classical group analysis.

Given a symmetry, one can
find integrating factors and first integrals, integrate the ODE in
quadratures, reduce ODE's order etc.
Infinitesimal symmetries turned out also to be a
useful tool for studying webs; see, for example,  \cite{MPPi},
where planar webs with infinitesimal symmetries were used for
construction of families of so-called exceptional webs.

As was mentioned above, the binary equation (\ref{binary}) defines a cubic ODE. Thus the obtained classification gives also a classification of some subclasses of cubic ODEs with a hexagonal webs of solutions. Studying of generic singular points of implicit ODEs was initiated
by Thom in \cite{Te}. For a generic quadratic ODE, normal forms were established by Davydov in \cite{Dn}. For cubic ODEs the classification problem
 becomes more complicated. It is clear that even a generic classification of cubic ODEs is not possible: the obstacle is the curvature.
  Moreover, Nakai showed that
the topological and analytic classifications are in fact the same in this case (see \cite{Nt}).
Even imposing the zero curvature condition will not compress the class of ODEs to guarantee a sensible classification. There is a partial result, which holds in both smooth and (real) analytic case. It is as follows. Our cubic ODE written as
$F(x,y,p)=0$ defines a surface $M$:
$$
M:=\{(x,y,p)\in {\mathbb K ^2 \times \mathbb P}^1(\mathbb K): \
F(x,y,p)=0\},
$$
where $(x,y,p)$ are coordinates in the jet space $J^1(\mathbb K,
\mathbb K)$ with $p=\frac{dy}{dx}$, $\mathbb K=\mathbb R$ or
$\mathbb C$.
The set of all points of the surface $M$, where the projection $\pi : M\to \mathbb K^2,\  (x,y,p) \mapsto
(x,y)$ is not a local diffeomorphism,
  is called a
 criminant.
Suppose the following regularity condition is imposed at
each point of the criminant:
$$
{\rm rank}((x,y,p)\mapsto (F,\frac{\partial F}{\partial p})=2.
$$
 This regularity condition
implies that the criminant and the surface $M$ are smooth. It turns out (see \cite{Ai}) that up to local diffeomorphism the above two examples
exhaust the list of normal forms. Namely these examples give normal forms if the projection $\pi$ has a cusp  at $(0,0)$.   To get normal forms for the case
of the fold singularity, where $p_1=p_2\ne p_3$, it suffices to pick up a
fold point on the discriminant curve of an example above and
rectify the integral curves corresponding to the root $p_3$. Thus a finite classification of implicit ODEs is possible if only one inpose some restriction. For an example of a class of implicit ODEs admitting such a slassification see  \cite{HIIY}. Of course, the obstacle of non-zero curvature does not make senseless the study of structurally stable properties of singular webs, such as the number of singular points, the sums of indices, etc. (see, for instance, \cite{Ys}, where the case of polynomial webs in $\mathbb {CP}^2$ was considered). The use of implicit ODEs as a tool for studying webs proved its efficiency also in studying abelian relations in non-singular case (see, for instance, \cite{Ha}).

The literature on the web geometry is immense. We mention here just a few references relevant to our case. In \cite{Nv} and \cite{Ne} the web structure was used for studying
geometric properties of differential equations. See \cite{Fc},\cite{Fi},\cite{Fh} for applications of 3-webs in mathematical physics, and  \cite{AGw}, \cite{PPi} for further references and surveys.

A note on the terminology used in this paper: as our considerations are local we will often omit the word "germ" in indicating the object under consideration.

\section{Chern connection and abelian relations}\label{Chern connection}
In this section we  present a formula for the Chern connection of
a 3-web, formed by solutions of an implicit cubic ODE. We use here
Blaschke's approach based on differential forms \cite{Be}.

\begin{definition}\label{def_hex}
Let $U\subset \mathbb C^2$ be an open set, where equation
(\ref{general_cub}) has 3 distinct roots $p_1,p_2,p_3$ and suppose
$U\ne \emptyset$. We say that this equation has a flat (or
hexagonal) 3-web of solutions if for
each point of $U$ there is a local
biholomorphism mapping the solutions of
(\ref{general_cub}) to three families of parallel lines.
\end{definition}

Let $p_1,p_2,p_3$ be the roots of (\ref{general_cub}) at a point
$(x,y)$ outside the discriminant curve.  1-forms vanishing on the
solutions can be chosen as follows
$$
\sigma_1=(p_2-p_3)(dy-p_1dx),\ \ \sigma_2=(p_3-p_1)(dy-p_2dx),\ \
\sigma_3=(p_1-p_2)(dy-p_3dx).\ \
$$
They are normalized to satisfy the condition
$$
\sigma_1+\sigma_2 + \sigma _3=0.
$$
Let us introduce an "area" form by  $$\Omega =\sigma_1\wedge
\sigma_2=\sigma_2\wedge \sigma_3=\sigma_3\wedge
\sigma_1=(p_1-p_2)(p_2-p_3)(p_3-p_1)dy\wedge dx.$$
 The Chern
connection form is defined as
$$
\gamma:
=h_2\sigma_1-h_1\sigma_2=h_3\sigma_2-h_2\sigma_3=h_1\sigma_3-h_3\sigma_1,
$$
where $h_i$ are determined by $$d\sigma_i=h_i\Omega.$$
The web is flat iff  the connection form
is closed: $d(\gamma)=0$. This implies $d\sigma_i=\gamma \wedge \sigma_i$. Defining
\begin{equation}\label{dk}
dk=-\gamma k,
\end{equation}  we introduce first integrals $u_i$ of the foliations  at least locally at a
regular point  by
\begin{equation}\label{du}
\begin{array}{c}
du_1=k\sigma_1=k(p_2-p_3)(dy-p_1dx),\\
du_2=k\sigma_2=k(p_3-p_1)(dy-p_2dx),\\
du_3=k\sigma_3=k(p_1-p_2)(dy-p_3dx).
\end{array}
\end{equation}
\noindent {\bf Remark.}  Let $\eta_1, \eta_2, \eta_3$ be germs of differential forms in $(\mathbb C^2,q_0)$ defining a flat 3 web and satisfying the following conditions:
\begin{itemize}
\item the forms are closed: $d(\eta_i)=0, \ \ i=1,2,3,$
\item the forms define the web: $\eta_i\wedge \sigma_i=0, \ \ i=1,2,3,$
\item the forms sum up to zero: $\eta_1 + \eta_2 + \eta_3=0,$
\end{itemize}
then these forms are proportional to $k\sigma_i$:  $\eta_i=c k\sigma_i, \ \ i=1,2,3$, $c=const$. One says that the space of abelian relations is one-dimensional for a hexagonal 3-web. I other words the first integrals summing up to zero are defined up to a constant factor. In what follows these first integrals are called {\bf abelian}.

\smallskip

To simplify
the final formulas we prefer to kill  the coefficient by $p^2$ in
equation (\ref{general_cub})  by a coordinate transform of the
form
\begin{equation}\label{chirn}
y=f(\tilde{x},\tilde{y}),\ x=\tilde{x},\ \ {\rm satisfying} \ \
3f_x(x,y)+a(x,y)=0.
\end{equation}
Therefore in what follows we often consider implicit ODEs without the
quadratic term:
\begin{equation}\label{CUB}
p^3+A(x,y)p+B(x,y)=0.
\end{equation}
By a coordinate transformation $y=F(X,Y), x=G(X,Y)$ the forms $\sigma_i$ are multiplied by the factor $\frac{(G_XF_Y-G_YF_X)^2}{G_X^3+aG_XG_Y^2-G_Y^3b}$,  where $a(X,Y)=A(G,F)$, and $b(X,Y)=B(G,F)$.

\begin{lemma}\label{changek}
Let $k(x,y)$ be a function not vanishing at $(0,0)$; then the following system of PDEs
\begin{equation}\label{normk}
\begin{array}{c}
k(G,F)(G_XF_Y-G_YF_X)^2=G_X^3+aG_XG_Y^2-G_Y^3b,\\
 \\
(3F_Y^2+aG_Y^2)F_X+2aF_YG_XG_Y+3bG_XG_Y^2=0
\end{array}
\end{equation}
 has a solution germ at  $(0,0)$ satisfying the conditions $(G_XF_Y-G_YF_X)\ne 0$, $F(0,0)=G(0,0)=0$.
\end{lemma}
{\it Proof:}  One easily checks the local solvability of the above system via the Cauchy-Kovalevskaya Theorem; locally the above system can be represented in Kovalevskaya form with respect to $F_X,G_X$ by adjusting Cauchy data. \hfill $\Box$\\

\bigskip

\noindent The following corollary is immediate.
\begin{lemma}\label{k1}
Suppose the Chern connection form is exact  $\gamma=d(f)$,
where the function $f$ is defined on some neighborhood  $U$ of a point on the discriminant curve. Then one can choose new local coordinates
to keep the coefficient by $p^2 $ to be zero and  simultaneously to ensure $k\equiv1$.
\end{lemma}
{\it Proof:} From (\ref{dk}) one has $k=exp(-f)\ne 0$.   Now choose $F,G$ to satisfy system (\ref{normk}) and $(G_XF_Y-G_YF_X)\ne 0$. The second equation of (\ref{normk}) ensures that the coefficient by $p^2$ remains zero. \hfill $\Box$\\
\smallskip

Computing the Chern connection form in terms of roots $p_i$ and using the Viete formulas one  gets
\begin{equation}\label{connection}
\gamma=\frac{(2A^2Ax-4A^2By+6ABAy+9BBx)}{4A^3+27B^2}dx+\frac{(4A^2Ay+6ABx+18BBy-9BAx)}{4A^3+27B^2}dy.
\end{equation}
Notice that this form can have a pole on the {\bf discriminant curve} $$\Delta:=\{(x,y): 4A^3(x,y)+27B^2(x,y)=0\}.$$
The condition $d(\gamma)=0$ gives the following
differential equation for the functions $A,B$ (\cite{Ai}):
\begin{equation}\label{PDE}
\begin{array}{l}
  (4A^3+27B^2)(9BA_{xx}-2A^2A_{xy}+6ABA_{yy}-6AB_{xx}-9BB_{xy}-4A^2B_{yy})+\\
  +108A^2BA_xB_y-108AB^2A_xA_y+162B^3A_y^2+40A^4A_yB_y-108A^2BA_x^2+\\
  +216A^2BB_y^2-36A^3B_xB_y+108A^2BA_yB_x-378AB^2A_yB_y-405B^2A_xB_x+\\
  -48A^3BA_y^2+8A^4A_xA_y+243B^2B_xB_y+84A^3A_xB_x+324ABB_x^2=0.\\
\end{array}
\end{equation}

\noindent {\bf Remark.} Computing normal forms, it is convenient to adapt local coordinates to the infinitesimal symmetry of the ODE. Then the equation can have a term with $p^2$ and a leading coefficient vanishing at the singular point. It is straightforward to derive the corresponding formulas for the connection form from  (\ref{connection}).

\section{Infinitesimal symmetries}\label{symsection}

Pick up a point $q_0$ on the discriminant curve and select some connected neighborhood $U$ of this point. At a point $q\in U\setminus \Delta$, equation (\ref{CUB}) implicitly  defines function germs $p_1,p_2,p_3$. Analytical continuation of these germs along all closed paths in $U$ passing through $q$ generates a subgroup of the group $S_3$ permuting the roots $p_i$. We call this subgroup a {\bf local monodromy group} of (\ref{CUB}) at $q_0$.

\smallskip

\noindent Notice that equation (\ref{CUB}) defines an  analytic set germ $\mathcal{A}$ in $(\mathbb C^5,0)$ by
\begin{equation}\label{setA}
p_1+p_2+p_3=0, \ \ p_1p_2+p_2p_3+p_3p_1=A(x,y), \ \ p_1p_2p_3=-B(x,y).
\end{equation}
We will need the following representation of functions holomorphic  on $\mathcal{A}$.
\begin{lemma}\label{preparationS3}
Suppose equation (\ref{CUB}) is irreducible over the ring of holomorphic function germs  $\mathcal{O}_0$ on $(\mathbb C^2,0)$  and the local monodromy group of (\ref{CUB}) acts on the roots as the permutation group $S_3$. Then each holomorphic function germ $F$ on the  analytic set germ $\mathcal{A}$
can be represented in the form
$$\textstyle
F=F_0(x,y)+p_1F_1(x,y)+p_2F_2(x,y)+p_1p_2F_3(x,y)+p_2^2F_4(x,y)+p_1p_2^2F_5(x,y),
$$
where $F_i$, $i=0,...,5$ are holomorphic function germs on $(\mathbb C^2,0)$. Moreover, this representation is unique.
\end{lemma}
{\it Proof:}  The existence of the representation
follows from Malgrange's Preparation Theorem. In fact, the identities
$$p_1^2=-p_1p_2-q_2^2-A,\ \ \ p_1^3=-p_1A-B,\ \ \ p_1^2p_2=-p_1p_2^2+B,\ \ \ p_2^3=-p_2A-B$$ imply
$\langle p_1,p_2 \rangle ^4 \subset \langle A,B \rangle$ and
$\mathcal{ O}_2(p_1,p_2)/ \langle A,B \rangle =\mathbb C
\{1,p_1,p_2,p_1p_2,p_2^2,p_1p_2^2 \}.$ To prove the uniqueness one applies all the permutations of $S_3$ to the representation of the zero function germ, normalize the results using the above identities and shows by direct computation that all $F_i$ are zero function germs. \hfill $\Box$\\

\smallskip

\noindent At each regular point one can choose any pair of the abelian first integrals defined by equations (\ref{du})  as local coordinates.
The symmetry algebra at this point in coordinates $u_1,u_2$ is generated by the following 3 vector fields
$$
\partial_{u_1},\ \ \partial_{u_2}, \ \ u_1\partial_{u_1}+u_2\partial_{u_2}.
$$
 If an operator $X$ is a symmetry then $X(u_i)=\varphi _i (u_i)$ for some function germs $\varphi _i$. As the space of abelian relations for a hexagonal 3-web is one-dimensional and the symmetry $X$ maps abelian relations into abelian relations, the functions $\varphi _i$ are linear:
 \begin{equation}\label{Xu}
 X(u_i)=Cu_i+c_i.
\end{equation}

\begin{lemma}\label{irreducibleS3}
Suppose $X$ is a symmetry of equation (\ref{CUB}) and the local monodromy group is $S_3$. Then $C\ne 0$ in the equality (\ref{Xu}).
\end{lemma}
{\it Proof:} Consider  a point $q_0=(x_0,y_0)$ not on the discriminant curve $\Delta$. Suppose $C=0$; then at least two of the constants $c_i$, say $c_1$ and $c_2$, do not vanish. Indeed, the corresponding first integrals are functionally independent at $q_0$ and a non-trivial symmetry operator cannot have 2 independent invariants. Equations (\ref{du}) imply
$$
c_1=k(p_2-p_3)(\eta-p_1\xi), \ \ \ c_2=k(p_3-p_1)(\eta-p_2\xi),
$$
where $X=\xi (x,y)\partial_x+\eta(x,y)\partial_y.$
Excluding the function $k$ gives $c_2(p_2-p_3)(\eta-p_1\xi)=c_1(p_3-p_1)(\eta-p_2\xi)$.
Rewriting this as $c_2\xi A+\eta (2c_1+c_2)p_1+\eta(2c_2+c_1)p_2-\xi(2c_1+c_2)p_1p_2-\xi(c_1-c_2)p_2^2=0$ and applying Lemma \ref{preparationS3} we get $X\equiv 0$.
\hfill $\Box$\\

\begin{theorem}\label{intrgrals via cub equation}
Suppose  ODE (\ref{CUB}) has a flat web of solutions and admits a symmetry operator $X$ such that $C\ne 0$ in the equality (\ref{Xu}).  Then one can choose  germs $I_i,\ i=1,2,3$ of first integrals of (\ref{CUB}) to satisfy $I_i=k^2 U_i$, where $dk=-k\gamma $,   $\gamma $ is the Chern connection form, and $U_i$ are the roots of the following cubic equation:
\begin{equation}\label{cub first integrals}
U^3-2\alpha U^2+\alpha ^2 U +\beta=0, \ \ {\rm where}
\end{equation}
 \begin{equation}\label{alphabeta}
 \alpha=\xi^2 A^2-3\eta^2 A-9\xi \eta B, \ \ \beta =(4A^3 + 27B^2)(\eta^3 + \xi^2\eta A+\xi^3 B)^2.
 \end{equation}
\end{theorem}
{\it Proof:} Let $q=(x_0,y_0)$ be a point not on the discriminant curve $\Delta$. Consider the germs of the first integrals $u_i$ at $q$ defined by (\ref{du}). Normalizing $X$ and adjusting integration constant in (\ref{du}) we have
$u_i=X(u_i)=k\sigma _i(X)$, i.e.
\begin{equation}\label{first integrals sym}
\begin{array}{c}
u_1=k(p_2-p_3)(\eta-p_1\xi ),\\
u_2=k(p_3-p_1)(\eta-p_2\xi), \\
u_3=k(p_1-p_2)(\eta-p_3\xi).
\end{array}
\end{equation}
Note that $I_i:=u_i^2$ are also first integrals. Using the Viete formulas for (\ref{CUB}) and relations (\ref{first integrals sym}) to calculate elementary symmetric function of $\frac{I_i}{k^2}$ one arrives at (\ref{cub first integrals}).
\hfill $\Box$\\

\bigskip

\noindent {\bf Remark.} Lie discovered that an explicit ODE in differentials $M(x,y)dx+N(x,y)dy=0$ with an infinitesimal symmetry $X$ has the integrating factor $\mu=\frac{1}{\xi M+\eta N }$, i.e. $d(\mu Mdx+\mu Ndy)=0$. The above Theorem gives an analog of this Lie result for implicit cubic ODEs.

\bigskip

\noindent {\bf Remark.}  Suppose the Chern connection form is exact  $\gamma=d(f)$,
where the function $f$ is defined on some neighborhood  $V$ of a point $q_0$ on the discriminant curve. Then one can normalize the web forms $\sigma _i$ to ensure $k=const$ (See Lemma \ref{k1}.) Now the first integrals can be chosen to satisfy equation  (\ref{cub first integrals}).

\bigskip

\noindent If $X$ is a symmetry of equation (\ref{CUB}) then the Lie derivative of the connection form $\gamma$  along the flow of $X$   vanishes. Therefore $\mathcal{L}_X(\gamma)=i_X(d(\gamma))+d(i_X(\gamma))=d(\gamma (X))=0$ since the connection form is closed. Thus  $\gamma (X)$ is constant:
 \begin{equation}\label{Xc}
 \gamma (X)=c.
 \end{equation}
\begin{theorem}\label{Xeuler} Suppose equation (\ref{CUB}) has a flat web of solutions,  $C\ne 0$ in the equality (\ref{Xu}), and a symmetry $X$ of (\ref{CUB})
vanishes at the
singular point $(0,0)\in \Delta$. Then the equation  is equivalent
to a weighted homogeneous ODE and the symmetry operator $X$ to an Euler vector field.
\end{theorem}
{\it Proof:} Choose a point $q=(x_0,y_0)$ not on the discriminant curve. Normalize the symmetry operator $X$ and the first integrals $u_i$ to satisfy $X(u_i)=u_i$. It is possible for $C\ne 0$.  Let us calculate the action of the symmetry operator $X$ on the functions $\alpha$ and $\beta$ defined by (\ref{alphabeta}):
$X(\alpha)=X(\frac{u_1^2+u_2^2+u_3^3}{k^2})=\frac{2u_1X(u_1)+2u_2X(u_2)+2u_3X(u_3)}{k^2}-2\frac{u_1^2+u_2^2+u_3^3}{k^3}X(k)=2(1-\frac{X(k)}{k})\alpha$. Since $\frac{X(k)}{k}=-(\gamma (X))=-c$ we have $X(\alpha)=2(1+c)\alpha$. Similarly $X(\beta)=6(1+c)\beta$.

We can suppose that the functions $\alpha$ and $\beta$ are functionally independent. If it is not true one can choose  a coordinate transform such that function germs  $\frac{u_i}{k}$ at $q$ are functionally independent. Therefore the functions $\alpha$, $\beta$  are also functionally independent. (This is a slightly modified version of Lemma \ref{changek}.)
Direct computation shows that this condition is equivalent to $c\ne -1$ in the formula (\ref{Xc}).

\smallskip

\noindent The operator $X$ vanishes at $(0,0)$ hence we can apply the following results of K.Saito (see \cite{Sq}).
\begin{enumerate}
  \item In suitable coordinates the operator $X$ can be written as a sum $X=X_s+X_n$ of an Euler operator $X_s$ (semi-simple in Saito's terminology) and a commuting with $X_s$ nilpotent operator $X_n=n_1(x,y)\partial _x +n_2(x,y)\partial _y)$ (i.e. all eigenvalues of the matrix $$\left( \begin{array}{cc}
                                                                                         \frac{\partial n _1}{\partial x} & \frac{\partial n_1}{\partial  y} \\
                                                                                         \frac{\partial n _2}{\partial x} & \frac{\partial n_2}{\partial  y}
                                                                                       \end{array}\right )
      $$ are zeros at $(0,0)$).
  \item Moreover,  the following two conditions are equivalent:

   a) $X(f)=\lambda f$,

   b) $X_s(f)=\lambda f$, $X_n(f)=0,$ \\where  $f$ is a function germ and $\lambda$ is a complex number.
\end{enumerate}
Thus we have from the condition b): $X_n(\alpha)=X_n(\beta)=0$. As the functions $\alpha$ and $\beta$ are functionally independent we get $X_n=0$. Hence the operator $X$ is an Euler operator in some new coordinates.
\hfill $\Box$\\

\begin{theorem} Suppose equation
(\ref{CUB}) has a flat web of solutions and admits  a symmetry $X$ at the
singular point $(0,0)\in \Delta$, the operator $X$ vanishes at this point, and the local monodromy group of (\ref{CUB}) is $S_3$. Then the equation  is equivalent
to a weighted homogeneous ODE and the symmetry operator $X$ to an Euler vector field.
\end{theorem}
{\it Proof:} The  claim follows from Lemma \ref{irreducibleS3} and Theorem \ref{Xeuler}.
\hfill $\Box$\\

\bigskip

\noindent {\bf Remark.} Unfortunately, the above Lemma is not true if the local monodromy group is smaller than $S_3$. It is not difficult to find  counter-examples:
\begin{itemize}
         \item equation $p^3=xy^6$ with the local monodromy group $Z_3$ admits the symmetry $y^2\partial _y$, which is not equivalent to an Euler vector field,
         \item equation $p(p^2-xy^4)=0$ with the local monodromy group $Z_2$ admits the symmetry $y^2\partial _y$, which is not equivalent to an Euler vector field,
         \item the web with abelian first integrals $u_1=\frac{1}{y-x^2}+x^2,\ u_2=\frac{1}{y-x^2},\   u_3=-\frac{2}{y-x^2}-x^2$ admits the symmetry $(y-x^2)^2\partial _y$, which is not equivalent to an Euler vector field. This web is defined by an ODE, which factors out into 3 linear in $p$ terms, i.e., its local monodromy group is trivial.
\end{itemize}

\begin{lemma}\label{first integrals for exact}
Let $q_0=(x_0,y_0)$ be a point on the discriminant curve $\Delta$.
Suppose the Chern connection form is exact  $\gamma=d(f)$,
where the function $f$ is defined on some connected neighborhood  $U$ of $q_0$. Then the abelian first integrals $u_i$ are integer algebraic over the ring of holomorphic function germs $\mathcal{O}_{q_0}$. They can be chosen to satisfy $u_i(q_0)=0$.
\end{lemma}
{\it Proof:}  Define $U_{\Delta}:=U\backslash \Delta$.
Then the connection form $\gamma$ is exact on $U_{\Delta}$.
Let $q\in U_{\Delta}$ be some point outside the discriminant curve
and $V$ a simply connected neighborhood of $q$ contained in $U_{\Delta}$, i.e.   $q\in V \subset U_{\Delta}.$
Select a path $\alpha: [0,1] \mapsto U$ connecting $q_0$ and $q$: $\alpha(0)=q_0$, $\alpha(1)=q$ and satisfying
$\alpha ((0,1])\in U_{\Delta}$.

Define functions $u_1,u_2,u_3: V \mapsto \mathbb C$ by equations (\ref{du}),
where  $k=\exp (-f)$ and $p_1,p_2,p_3: V \mapsto  \mathbb C$ are functions implicitly defined by equation (\ref{CUB}).
Then $u_1,u_2,u_3$ are well-defined up to a choice of the initial values $u_1(q),\ u_2(q),\ u_3(q)$.
Let us fix them by
$$
\begin{array}{c}
u_1(q)=\int_{\alpha}k(p_2-p_3)(dy-p_1dx),\\
 u_2(q)=\int_{\alpha}k(p_3-p_1)(dy-p_2dx),\\
  u_3(q)=\int_{\alpha}k(p_1-p_2)(dy-p_3dx).
\end{array}
$$
The analytical continuation of $u_i$ along all the paths contained in $U_{\Delta}$ gives multivalued functions $\tilde{u}_i$ on $U_{\Delta}$. Due to the choice of initial conditions one has
$$
\tilde{u}_1+\tilde{u}_2+\tilde{u}_3=0.
$$
Moreover, these initial conditions also imply that the functions
$$
f:=\tilde{u}^2_1+\tilde{u}^2_2+\tilde{u}^2_3, \ \ \ h:=\tilde{u}^2_1\tilde{u}^2_2\tilde{u}^2_3
$$
are one-valued on $U_{\Delta}$. In fact, the analytic continuation along each closed path in $U_{\Delta}$
induces a permutation of roots $p_1,p_2,p_3$. This permutation generates an action on the differentials $d(u_1),d(u_2),d(u_3)$. For example, to the cycle $(2,3,1)$ there corresponds the same permutation of the differentials, while the permutation $(2,1,3)$ generates the following transformation: $$d(u_1)\to -d(u_2),\ d(u_2)\to -d(u_1),\ d(u_3)\to -d(u_3).$$  On the other hand, due to the choice of
the initial values of $u_i$ the action on $\tilde{u}_1,\tilde{u}_2,\tilde{u}_3$ coincides with the action on the differentials. Moreover, being bounded, the functions $f,h$ are holomorphic on the whole neighborhood $U$ by Riemann theorem.
Therefore each of the functions $\tilde{u}_i$ is integer over the ring $\mathcal{O} (U)$ of functions analytical on $U$ as satisfying the following equation
$$
u^6-fu^4+\frac{f^2}{4}u^2-h=0
$$
Differentiating the function $f$ and using (\ref{du}) one shows that the  functions $\tilde{u}_1,\tilde{u}_2,\tilde{u}_3$ are well
defined meromorphic functions on the germ of analytic set $\mathcal{A}$  determined by equations (\ref{setA}).
Further, being integer also over $\mathcal{O}(A)$ these functions are in fact holomorphic on $\mathcal{O}(A)$. \hfill $\Box$\\

\bigskip

\noindent
According to the classical Lie results the components of a
symmetry operator $X$ satisfy  a system of linear PDEs. In a
neighborhood of a regular point the space of solutions to this
system is 3-dimensional.  When there exists a solution that can be
extended to a neighborhood of a point on the discriminant curve $\Delta$?
A sufficient condition  gives the following theorem.
\begin{theorem} Let $q_0=(x_0,y_0)$ be a point on the discriminant curve $\Delta$.
 Suppose the Chern connection form is exact  $\gamma=d(f)$,
where the function $f$ is defined on some neighborhood  $U$ of $q_0$. Then the dimension of the symmetry algebra of  equation (\ref{CUB}) at $q_0$ is
\begin{itemize}
\item at least 1, if all 3 roots coincide and the local monodromy group is $S_3$,
\item at least 2, if only 2 roots coincide and the local monodromy group is $Z_2$.
\end{itemize}
\end{theorem}
{\it Proof:} Define the first integrals as in Lemma \ref{first integrals for exact}.

{\it $\bullet$ Triple root.}
Each holomorphic function germ on $\mathcal{A}$ can be written in the normal form given by  Lemma \ref{preparationS3}.
Using the symmetry properties of $\tilde{u}_1,\tilde{u}_2,\tilde{u}_3$ under the permutations of the roots one gets
$$
\tilde{u}_1=(p_2-p_3)(M(x,y)-p_1L(x,y)),\ \ \ \tilde{u}_2=(p_3-p_1)(M(x,y)-p_2L(x,y)),
$$
where $M$ and $L$ are aholomorphic on $U$ (compare with \cite{Ai}).
Define
\begin{equation}\label{define symmetry X}
\xi =\frac{1}{k(p_1-p_2)}\left( \frac{\tilde{u}_2}{p_3-p_1}-\frac{\tilde{u}_1}{p_2-p_3} \right ), \ \ \ \eta=\frac{1}{k(p_1-p_2)}\left( \frac{p_1\tilde{u}_2}{p_3-p_1}-\frac{p_2\tilde{u}_1}{p_2-p_3} \right ).
\end{equation}
It is immediate that the functions $\xi=\frac{L}{k}$ and $\eta=\frac{M}{k}$ are well defined on $U$. Therefore the vector field
$X=\xi(x,y)\partial_x +\eta (x,y) \partial _y$ is a symmetry of our ODE, as its action on the first integrals satisfies $X(u_1)=u_1$, $X(u_2)=u_2$ due to the equalities (\ref{define symmetry X}) and
(\ref{du}).
\smallskip

{\it $\bullet$ Double root.} Suppose $p_1,p_2$ satisfy an irredicible quadratic equation at $0$ and  $p_1=p_2\ne p_3$. Then $a:=p_3$ is a holomorphic function germ on $(U,q_0)$ and $a(q_0)\ne 0$ since $p_1+p_2+p_3=0$. The function germ $2p_3^2 +p_1p_2=2a^2+p_1p_2$ is also holomorphic.
Moreover, it does not vanish at $q_0$: $2a^2+p_1p_2|_{q_0}=\frac{9}{4}a^2(q_0)$. Then the vector field
$$
X_1=\frac{\partial_x+p_3\partial_y}{k(2p_3^2 +p_1p_2)}
$$
is an infinitesimal symmetry. Indeed, its action on the first integrals is the following: $X_1(u_1)=-1$, $X_1(u_2)=1$.
The second symmetry $X_2$ is defined by the same formula as for the case of triple root. To check that the vector field (\ref{define symmetry X}) is well defined on $(U,q_0)$ write
$$
\tilde{u_1}(p)=R(x,y)+p_1S(x,y)
$$ instead of the normal form given by Lemma
\ref{preparationS3}, observe that $\tilde{u}_2(p)=aSR-R+p_1S$ due to the permutation symmetry properties, and substitute these expressions into (\ref{define symmetry X}).
One immediately checks that this vector field satisfies  $X_2(u_1)=u_1$, $X_2(u_2)=u_2$. On some neighborhood $V\subset U$ of a point $q\ne q_0$,  one can rewrite the symmetry operators as
$X_1=\partial_{u_2}-\partial_{u_1}$, $X_2=u_1\partial_{u_1}+u_2\partial_{u_2}$, i.e. they are linearly independent.
\hfill $\Box$\\
 \\

\noindent {\bf Remark.}
Note that as $U_{\Delta}$ in Lemma \ref{first integrals for exact} is not simply connected the Poincar\'e lemma is not applicable.
Moreover, in the above proof we need the function $f$ to be defined also on $\Delta$, not only on $U_{\Delta}$.\\

\noindent {\bf Remark.} Note that in the above proof we essentially use the relations on $F_i$ (case $S_3$) or $R,S$ (case $Z_2$) derived from the symmetry properties. It seems that the regularity of the Chern connection $\gamma$ does not suffice for existence of symmetry. Consider, for example, the web with the abelian integrals $u_1=y$, $u_2=y+y^2+yx^2$, $u_3=-u_1-u_2=-2y-y^2-yx^2$. Then $\gamma \equiv 0$ but this web does not admit any symmetry at $(0,0)$. In fact, this point is singular for the discriminant curve, therefore $X|_0=0$. Hence $X(u_i)=cu_i$. Substituting $X=\xi \partial _x +y \partial _y$ one arrives at  $y+2x\xi=0$.

\section{Weighted homogeneous ODEs}
The simplest case of a
symmetry with an isolated fixed point is a scaling. This case is
also the most interesting for applications in physics: to define a
structure of a Frobenius 3-fold, solutions of associativity
equations (\ref{ass}) or (\ref{ass2}
) must be weighted homogeneous.
\begin{definition}
We call an implicit ODE weighted homogeneous if its web of
solutions is invariant with respect to the flow of non-trivial
Euler vector field (\ref{euler}). The numbers $w_x,w_y$ are called
weights.
\end{definition}
In this section we present a classification of singularities of weighted homogeneous implicit ODEs with hexagonal 3-webs of solutions. The equivalence
relation used is the group of biholomorphism germs preserving the
origin.\\
 \noindent We call the Euler vector field
{\bf parabolic, hyperbolic} or {\bf elliptic} if $w_xw_y=0$,
$\frac{w_x}{w_y}<0$ or $\frac{w_x}{w_y}>0$ respectively. Later we will see that the weights  $w_x, y_y$ can be chosen to be rational therefore the above definition makes sense.

If exactly two web directions coincide and distinct 2 directions are those of the coordinate axes then one has $K_3(0,0)=K_0(0,0)=0$ for  equation (\ref{binary}). At least one of the values $K_2(0,0)$ or $K_1(0,0)$ does not vanish. Without loss of generality we assume $K_2(0,0)\ne 0$ and rewrite (\ref{binary}) as
$$
f(x,y)p^3+p^2+g(x,y)p+h(x,y)=0,
$$
where $f(0,0)=h(0,0)=0$. There exists a coordinate transform $y=\varphi (\bar{x},\bar{y}),\ x=\bar{x}$ respecting the symmetry and "killing" the coefficient $g$. Indeed, there is a function germ $\varphi$ satisfying the following PDE and Cauchy data:
$$
g(\bar{x},\varphi)+2\varphi_{\bar{x}}+3f(\bar{x},\varphi)\varphi^2_{\bar{x}}=0, \ \ \ \varphi(0,\bar{y})=\bar{y}.
$$
The above problem is symmetric with respect to the symmetry
$E+w_y\varphi \partial_{\varphi}$ and so is its unique solution.
Interchanging again the coordinates we arrive at
\begin{equation}\label{nonmonicCUB}
F(x,y)p^3+p+G(x,y)=0
\end{equation}
with $F(0,0)=G(0,0)=0$.

If one of the coordinate axes is different from the foliation directions (for example, if all 3 directions coincide at $(0,0)$) then this coordinate axis is transverse to all foliations. Interchanging if necessary the exes we can reduce ODE (\ref{binary}) to the monic form (\ref{CUB}). Using the above arguments one checks that the coordinate transformation (\ref{chirn}) can
be chosen to respect the weighted homogeneity.

Our approach to classification is the following. First we apply a coordinate transform respecting the scaling symmetry to simplify a non-zero coefficient.
To ensure the existence of the transform the following Lemma will be used.
\begin{lemma}\label{productLM}
Suppose $f:(\mathbb C,0)\to (\mathbb C,1)$ is holomorphic and $n,m\in \mathbb N$. Then the ODE $$\psi ^nf(\psi)\left(\frac{d\psi}{dt}\right)^m=t^n$$ has an holomorphic solution $\psi: (\mathbb C,0) \to (\mathbb C,0)$ with $\left.\frac{d\psi}{dt}\right|_{t=0}\ne 0$.
\end{lemma}
{\it Proof:} Integrating  $\psi ^{\frac{n}{m}}f(\psi)^{\frac{1}{m}}d\psi =t^{\frac{n}{m}}dt$ we get $\psi ^{\frac{n}{m}+1}\tilde{f}(\psi) =t^{\frac{n}{m}+1}$, where $\tilde{f}(0)\ne 0$. Therefore the equation $\psi \tilde{f}(\psi)^{\frac{m}{n+m}} =t$ gives the desired solution. \hfill $\Box$\\

\smallskip
\noindent The second step is to apply a
 monomial substitutions
\begin{equation}\label{equivalence}
x=m\bar{x}^{\alpha},\ \ \ y=l\bar{y}^{\beta}
\end{equation}
with constant $m,l$ and suitable rational $\alpha, \beta$  (in
particular, scalings $x=m\bar{x},\ \ \ y=l\bar{y}$). Observe that these substitutions obviously
preserve the class of weighted homogeneous ODEs.

Finally, we analyze the condition $d(\gamma)=0$ and  present the corresponding normal form applying the inverse of (\ref{equivalence}).
 We also use these substitutions as a tool
reducing lists of normal forms.

\subsection{Parabolic case $w_xw_y=0$}

In this subsection we give normal forms for the case when one of the weights vanishes.

\begin{theorem}\label{parabolicTH}
Suppose one of the weights of the weighted homogeneous ODE (\ref{binary}) with a flat web of
solutions vanishes, then for some nonnegative integer $N$ and constant $L_0,L_1$ the equation is equivalent to one from the following list:
$$
\begin{array}{cl}
1) & p^3=x^{N}y^3\\
2) & p^3+x^Ny^2p=\frac{2}{\sqrt{27}}x^{\frac{3}{2}N}y^3\tan(L_0x^{1+\frac{N}{2}}+L_1)\\
3) & p^3+x^{N+1}y^2p+\frac{2x^{\frac{3}{2}(N+1)}y^3}{\sqrt{27}\tan\left(L_0x^{1+\frac{N+1}{2}}\right)}=0\\
4) & p(x^2y^Np^2+1)=0.
\end{array}
$$
For the form 2)  $L_1\ne \frac{\pi}{2}$, $L_1=0$ for odd $N$, and with $L_1\ne 0$ one can choose $0 \leq \arg(L_1) < \pi$.
\end{theorem}
{\it Proof:} { $\bullet$ \it Monic case}. It is easy to see that for $w_y=0$, $w_x\ne
0$ there is no smooth monic cubic  weighted homogeneous ODE: the
coefficients must have the forms  $A=a(y)/x^2$, $B=b(y)/x^3$.

\noindent  Let us normalize the weight:
$w_y=1$. Then $w_p=1$ and equation (\ref{CUB}) becomes
\begin{equation}\label{pCUBg}
p^3+a(x)y^2p+y^3b(x)=0.
\end{equation}
Let $a(x)\not \equiv 0$ then $a(x)=x^k\alpha(x)$ with $\alpha(0)\ne 0$ for
some nonnegative integer $k$. Consider the following coordinate
transform:
\begin{equation}\label{parabolicCoordinateTransform}
y=\tilde{y},\ \ \  x=\varphi (\tilde{x}).
\end{equation}
Note that it preserves the Euler vector field $E=y\partial_y$. In
the new coordinates the ODE takes the form
$$
\tilde{p}^3 +\alpha (\varphi)\varphi ^k (\varphi
')^2\tilde{y}^2\tilde{p}+b(\varphi)(\varphi ')^3y^3=0.
$$
Choose the function $\varphi (\tilde{x})$ to satisfy $\alpha
(\varphi)\varphi ^k (\varphi ')^2=\tilde{x}^k$ and $\varphi
(0)=0$ (Lemma \ref{productLM}), then (\ref{parabolicCoordinateTransform}) with this
$\varphi$ correctly defines a coordinate transform bringing
(\ref{pCUBg}) to the form
$$
p^3+x^ky^2p+y^3b(x)=0.
$$ This equation has a flat web of solutions, which
impose a second order ODE on $b(x)$. To simplify the analysis let
us "kill" the factor $x^k$ by $p$ via a suitable substitution (\ref{equivalence}). Then one arrives at
$$
p^3+y^2p+y^3\beta (x)=0
$$
with the connection form
$$
\gamma =\frac{9\beta '(x)d(x)}{27\beta^2(x)+4}+\frac{(6\beta
'(x)+54\beta^2(x)+8)dy}{(27\beta^2(x)+4)y}.
$$
This form is closed, which implies
$$
\frac{3\beta '(x)+27\beta^2(x)+4}{27\beta^2(x)+4}=const.
$$
Integration gives $p^3+y^2p=\frac{2}{\sqrt{27}}y^3\tan(c_0x+c_1)$.
Now applying (\ref{equivalence}) one gets the forms 2) and 3). The substitution $x\to \alpha x$ with $\alpha ^{1+\frac{N}{2}}=-1$ changes $L_1$ for $-L_1$.

\noindent If $a(x)\equiv 0$ then $b(x)=x^k\beta(x)$ with $\beta(0)\ne
0$ for some nonnegative integer $k$. Then coordinate
transformation (\ref{parabolicCoordinateTransform}) reduces the
equation to
$$
p^3=x^ky^3.
$$

\noindent{$\bullet$ \it Non-monic case.} If $w_x=0$ normalize the weight $w_y$ as above.
Then $w_p=1$ and the functions $F,G$ in (\ref{nonmonicCUB}) take the forms $F(x,y)=\frac{f(x)}{y^2}$, $G(x,y)=yg(x)$.
Thus $f\equiv 0$ and two directions coincide identically.

Hence the symmetry is $-x\partial _x+p\partial_p$. This implies $F(x,y)=x^2f(y)$ and $G(x,y)=\frac{g(y)}{x}\equiv 0$. It is immediate that the web of solutions of equation
$$
p(x^2f(y)p^2+1)=0
$$
is hexagonal: the first integrals can be chosen as follows $I_1=u(y)+ \ln x$, $I_2=u(y)-\ln x$, $I_3=-2u(y)$ for some function $u$ determined by the quadratic factor of the ODE.
Thus $I_1+I_2+I_3=0$.

For some function $\alpha$ with $\alpha (0)\ne 0$ holds true $f(y)=y^k\alpha (y)$.
 Choosing function $\varphi$ in the coordinate transform
$
x=\tilde{x},\ \ \  y=\varphi (\tilde{y})
$
to satisfy $\varphi^k\alpha (\varphi)(\varphi ')^2=y^k$ one arrives at $p(x^2y^kp^2+1)=0$.
\hfill $\Box$\\

\subsection{Elliptic case $w_x/w_y>0$}
\begin{theorem}\label{ellipticTH}
If the weights of the weighted homogeneous ODE (\ref{binary}) with a flat web of
solutions satisfy $w_1w_2>0$, then the equation is obtained from one of
the following list via suitable substitution (\ref{equivalence}).
$$
\begin{array}{cl}
1) & p(p-1)(p+1)=0 \\
2) & p^3+xp-y=0  \\
3) & p^3+2xp+y=0\\
4) & p^3+yp+\frac{1}{9}xy=0\\
5) & p^3+yp-\frac{2}{9}xy=0\\
6) & p^3-x^2p+x(\frac{2}{\sqrt{3}}x^2-\frac{8}{3}y)=0\\
7) & (p+\frac{x}{3})(p^2-\frac{1}{3}xp-y+\frac{x^2}{9})=0\\
8) & (p-\frac{2}{3}x)(p^2+\frac{2}{3}xp+y-\frac{2}{9}x^2)=0\\
9) & (p-\frac{2}{3}x)(p^2+\frac{2}{3}xp-2y-\frac{2}{9}x^2)=0\\
10) & p^3+(6y+3x^2)p+2x^3=0\\
11) & p^3-\frac{3}{2}(y+x^2)p-x(y+x^2)=0\\
12) & p^3-3(y+\frac{x^2}{2})p-x(y+\frac{x^2}{2})=0\\
13) & p^3+(3y-\frac{3}{2}x^2)p+x(x^2-3y)=0\\
14) & p^3-xyp+y(y+\frac{x^3}{27})=0\\
15) & p^3+\frac{1}{2}xyp+y(y+\frac{x^3}{27})=0\\
16) & p^3+2xyp+y(y-\frac{8}{27}x^3)=0\\
17) & p^3+x(y-\frac{2}{9}x^3)p+(y-\frac{2}{9}x^3)^2=0\\
18) & p^3+\frac{5}{2}x(y
 -\frac{25}{18}x^3)p+(y-\frac{25}{18}x^3)^2=0\\
19) & p^3+x(3y - x^3)p+(y-\frac{x^3}{3})(y-\frac{4x^3}{3})=0\\
20) &
 p^3+x(\frac{3}{2}y-\frac{2}{25}x^3)p+(y-\frac{2}{15}x^3)(y-\frac{4}{75}x^3)=0\\
 21)& p^3-x(2y +
 \frac{x^3}{9})p+y^2+\frac{x^6}{81}+\frac{4}{9}yx^3=0\\
 22) &
 p^3+x(\frac{5}{2}y-\frac{2}{9}x^3)p+y^2+\frac{4}{81}x^6-\frac{7}{9}yx^3=0\\
 23) &
 p^3+4x(y-\frac{4}{9}x^3)p+y^2+\frac{64}{81}x^6-\frac{32}{9}yx^3=0\\
24) & p^3+x(y+\frac{x^3}{36})p+y^2+\frac{x^6}{324}+\frac{yx^3}{18}=0\\
 25) & p^3-x(\frac{y}{2}-\frac{5}{2^8 3^2}x^3)p+y^2-\frac{x^6}{2^{11}3^4}-\frac{y x^3}{2^5
 3^2}=0\\
 26)& p^3+x(\frac{5}{2}y + \frac{5^3}{2^8 3^2}x^3)p+y^2-\frac{5^4 13}{2^{11} 3^4}x^6-\frac{5^2 7}{2^5
 3^2}yx^3=0.
\end{array}
$$

\end{theorem}
{\it Proof:} { $\bullet$ \it Monic case}. Let
us normalize the weights $w_x,w_y$ to satisfy $w_p:=w_y-w_x=1$. Then the functions
$A$ and $B$ have the weights $2$ and $3$ respectively. Let their
Tailor series be
\begin{equation}\label{TailorAB}
A(x,y)=\sum a_{k,l}x^ky^l,\ \ B(x,y)=\sum b_{m,n}x^my^n.
\end{equation} From the homogeneity condition one has :
$$
a_{k,l}[2-(kw_x+lw_y)]=0, \ \ b_{m,n}[3-(mw_x+nw_y)]=0,\  {\rm
where} \ k,l,m,n\in \mathbb N_0.
$$
This implies
\begin{equation}\label{pw}
w_x(k+l)+l=2, \ \ w_x(m+n)+n=3,
\end{equation}
 for non-vanishing coefficients
$a_{k,l}$, $b_{m,n}$. Moreover, the weights are rational. Now we can rewrite the condition  $w_x/w_y>0$ as $w_x(w_x+1)>0$. If $w_x<-1$ then $w_y<0$ and the weight of $B$ cannot be equal to 2: the weight of all the terms of $B$ are negative. Therefore $w_x>0$ and $l=0,1$ and $n=0,1,2.$ Solving (\ref{pw}) one gets the following pre-normal forms
$$
A(x,y)=ax^{2N-1}, \ \ \ B(x,y)=byx^{2(N-1)} \ \ \ \ {\rm or}
$$
$$
A(x,y)=ax^{2N}+byx^{N-1}, \ \ \
B(x,y)=rx^{3N}+syx^{2N-1}+ty^2x^{N-2},
$$
where $a,b,r,s,t$ are constant, $N\in \mathbb N$, and $t=0$ for
$N=1$. For $w_y=w_x$ the normalization
$w_p=1$ is not possible. In this case the coefficients $A$ and $B$
are constant  and the web is regular. Substituting the pre-normal forms into (\ref{PDE}) and collecting
similar terms one arrives at a set of polynomial equations for the
coefficients of the pre-normal forms. Fortunately, these equations
can be explicitly solved. Normalization of the obtained
coefficients by means of substitutions (\ref{equivalence}) gives a
finite list of normal forms, each of them giving an infinite
series of ODEs after applying (\ref{equivalence}) with suitable
parameters.\\

{$\bullet$ \it Non-monic case.} Here it is convenient to normalize the weights to satisfy $w_p=-1$. As above one shows that the weights are rational. From $w_x=w_y+1$ it follows that $w_y(1+w_y)>0$. If $w_y<-1$ then the weight of $F$ in (\ref{nonmonicCUB}) is 2. Since the weights of all monomials are negative, one has $F\equiv 0$ and two directions coincide identically. Therefore $w_y>0$, which implies $G\equiv 0$ as having the weight $-1$. Now our equation takes the form
$
p(F(x,y)p^2+1)=0,
$
where the weight of $F$ is 2. Expanding $F$ into Taylor series
$F(x,y)=\sum f_{k,l}x^ky^l,$
one concludes that $(k+l)(w_y+1)=2-k$ and $k\le 1$. Hence $F=f_0y^{l_1}+f_1xy^{l_0}$ for some natural $l_1$ and nonnegative integer $l_0$.

If $f_1=0$  then a suitable substitution (\ref{equivalence}) brings the equation to form 1). If $f_1\ne 0$  substitution (\ref{equivalence}) permits to transform the equation to
$p((cy^n+x)p^2+1)=0$ with a natural $n$. The condition $d\gamma =0$ implies $c=0$ and we get again the form 1) after substitution  (\ref{equivalence}) with
 $\alpha = \frac{1}{2}$.
\hfill $\Box$\\
 \smallskip

\noindent {\bf Remark.}
Each equation of the  list from Theorem \ref{ellipticTH} generates an infinite series of
weighted homogeneous equations with a flat web of solutions.

\subsection{Hyperbolic case $w_x/w_y<0$}
For analysis of normal forms, it is convenient to introduce an invariant variable of the form
$s:=x^{n}y^{m}$, where $n,m\in \mathbb N$ will be defined below to satisfy $E(s)=0$. To simplify the ODE we use the following coordinate transformation
\begin{equation}\label{htransform}
y=YQ(S),\ \ \ x=XR(S),
\end{equation}
where $S=X^nY^m$, $Q(0)\ne 0$ and $R(0)\ne 0$. One easily
checks that this transform is invertible as the relation
$s=SR^n(S)Q^m(S)$ allows one to find locally $S$ as a function
of $s$. Moreover, this transform preserves the Euler vector field.
The substitution (\ref{htransform}) transforms an "hyperbolic" monic ODE
$$
p^3+\frac{y}{x}\sigma(s)p^2+\left(\frac{y}{x}\right)^2\alpha (s)p+\left(\frac{y}{x}\right)^3\beta(s)=0,
$$
where $\sigma ,\alpha , \beta $ are analytic with zeros of suitable orders, to
$$
P^3+\frac{Y}{X}\tilde{\sigma}(S)P^2+\left(\frac{Y}{X}\right)^2\tilde{\alpha} (S)p+\left(\frac{Y}{X}\right)^3\tilde{\beta}(S)=0.
$$
The equivalence of 2 equations amounts to 3 ODEs for 2 functions $Q,R$. In the non-monic case we also have 3 ODEs for $P,Q$.
To eliminate redundant normal forms forms, we use the following Lemma.
\begin{lemma}\label{hreduce}
Suppose one can locally find the derivatives of $Q,R$ as analytic functions of $P,Q$ and $s$ at a point $(P_0,Q_0,0)$ with $P_0\ne 0$, $Q_0\ne 0$, from the system of the above 3 equations, where these equations are satisfied at $(P_0,Q_0,0)$. Then the corresponding ODEs are equivalent.
\end{lemma}
{\it Proof:} The system of those 3 equations is compatible, the web flatness being the compatibility condition. Therefore the equations are satisfied identically with the found $Q',R'$, being satisfied at one point. \hfill $\Box$\\
\begin{theorem}\label{hyperbolicTH}
If the weights of the weighted homogeneous ODE (\ref{binary}) with a flat web of
solutions satisfy $w_1w_2<0$, then for some non-negative integers $n_0,m_0,l_0$ and constant $L\ne \frac{1}{3(1-2k)}$, $k\in \mathbb N$, the equation is equivalent to one of the
following list:
$$
 \begin{array}{ll}
1)&p^3+x^{n_0}y^{3+m_0}p+x^{\frac{3n_0}{2}}y^{\frac{9+3m_0}{2}}U\left(\left[\frac{2(m_0+1)}{n_0+2}\right]x^{1+\frac{n_0}{2}}y^{\frac{1+m_0}{2}}\right)=0, \\
2)&p^3+x^{n_0+1}y^{3+m_0}p-\frac{x^{\frac{3(n_0+1)}{2}}y^{\frac{9+3m_0}{2}}}{V\left(\left[\frac{2(m_0+1)}{n_0+3}\right]x^{1+\frac{n_0+1}{2}}y^{\frac{1+m_0}{2}}\right)}=0, \\
3)&p(x^{3+n_0}y^{m_0}p^2+1)=0,\\
5)&\left[\frac{1+m_0}{1+n_0}\right]^{5+l_0}x^{N}y^{M}W\left(\left[\frac{1+m_0}{1+n_0}\right]x^{1+n_0}y^{1+m_0}\right)p^3+p+x^{n_0}y^{2+m_0}=0,
\end{array}
$$
where $N=(1+n_0)(3+l_0)+2$, $M=1+l_0+(l_0+3)m_0$.

\noindent The function $U(T)$  either vanishes identically or is defined, with suitable constants $C_1,C_2$, by the relations
\begin{equation}\label{UPQ}
\begin{array}{ll}
\frac{T}{3\sqrt3L}=\frac{f'\left(-\arctan\left(\frac{3\sqrt3}{2}U\right)\right)}{f\left(-\arctan\left(\frac{3\sqrt3}{2}U\right)\right)},& f(z)=\cos^{-\mu} (z)[C_1P^{\mu}_{\nu}(\sin z)+ C_2Q^{\mu}_{\nu}(\sin z)]\\
f\left(-\arctan\left(\frac{3\sqrt3}{2}U(0)\right)\right)\ne 0, & f'\left(-\arctan\left(\frac{3\sqrt3}{2}U(0)\right)\right)=0.\\
\end{array}
\end{equation}
The initial value of $U$ vanishes $U(0)=0$ if at least one of the numbers $n_0$, $m_0+1$ is odd. If $U(0)\ne 0$ then one can choose $0 \leq \arg(U(0)) < \pi$.

\noindent The functions  $V(T)$ is defined by the relations
\begin{equation}\label{VPQ}
\begin{array}{lll}
\frac{T}{3\sqrt3L}=\frac{f'\left(-\arctan\left(\frac{2}{3\sqrt3}V\right)\right)}{f\left(-\arctan\left(\frac{2}{3\sqrt3}V\right)\right)},& f(z)=\sin^{\mu} (z)P^{\mu}_{\nu}(\cos z).& \\
\end{array}
\end{equation}
If $L= -\frac{2}{3}$, the numbers $n_0,m_0$ are  odd, and $n_0\ge3$, the function  $V$ is either as in (\ref{VPQ}) or is defined by the relations
\begin{equation}\label{VPQ2}
\begin{array}{ll}
\frac{T}{3\sqrt3L}=\frac{f'\left(-\arctan\left(\frac{2}{3\sqrt3}V\right)\right)}{f\left(-\arctan\left(\frac{2}{3\sqrt3}V\right)\right)},& f(z)=\sin^{\mu} (z)[P^{\mu}_{\nu}(\cos z)+ Q^{\mu}_{\nu}(\cos z)].\\
\end{array}
\end{equation}
In the above formulas, $P^{\mu}_{\nu}(z),Q^{\mu}_{\nu}(z)$ are the Legendre functions for $\mu=\frac{1}{2}(1-\frac{1}{3L})$, $\nu=\frac{1}{2}(\frac{1}{L}-1)$.

\noindent Finally, $W(t)$ with $W(0)\ne 0$  is a solution to
$$[2(1-s)+3s^{4+l_0}(3s-1)W]\frac{dW}{ds}=\frac{5+l_0}{2}[3s^{3+l_0}(2+3s)W+4]W.$$
\end{theorem}
{\it Proof:} Let us normalize the weights to satisfy $w_y>0$, $w_x<0$,
$w_p:=w_y-w_x=1$.

\noindent {$\bullet$ \it Monic case}.
 If there is at least one term in the
Tailor expansion of $A(x,y)$ (or $B(x,y)$), say $a_{ij}x^jy^j$,
then $iw_x+jwy=2$ (or $iw_x+jwy=3$). With $w_x=w_y-1$ one
concludes immediately that the normalized weights are rational.
Moreover, choosing coprime numbers $q,r\in \mathbb N$  such that $w_y=\frac{q}{r}$ one has $q<r$.

 Let us introduce an invariant variable
$s:=x^qy^{r-q}$, where $r,q$ are defined above. Consider the Tailor
expansions (\ref{TailorAB}) of the coefficients $A,B$. Then the
weight of $A$ is equal to 2, i.e., $E(A)=2A$ and the weight of $B$
is equal to 3: $E(B)=3B$. This allows one to represent $A,B$ in
the following form:
$$
A(x,y)=x^{ql-2}y^{2+l(r-q)}a(s),\ \ \ B(x,y)=x^{qt-3}y^{3+t(r-q)}b(s),
$$
where $l,t\in \mathbb N$ satisfying $lq\ge 2$, $tq\ge 3$. First consider the case when $a$ does not
vanishes identically.\\

{\it Case $a(s)\not \equiv 0$.} As $a(s)$ does not vanish identically, we can rewrite the equation
in the form
$$
p^3+\left(\frac{y}{x}\right)^2s^k\alpha (s)p+\left(\frac{y}{x}\right)^3\beta(s)=0
$$  with $\alpha (0)\ne 0$.
One can choose the functions $Q,R$ of the transform (\ref{htransform}) to keep the
coefficient by $p$ to be zero and to make  $\alpha(s)$ constant.
This amounts to a cumbersome but direct verification that a system of two ODEs for
$Q,R$ locally has a suitable solution.(Compare with the parabolic case, where we have only one
¨free¨ function $\varphi$ and one ODE for it.) Thus, we can assume
that $\alpha (s)\equiv 1$. Further, applying substitution
(\ref{equivalence}) we kill the dependence on $x$ of
$A(x,y)$ and normalize it to $y^3$. Now $s=x^2y$ and our equation
is
$$
p^3+y^3p+xy^5b(x^2y)=0.
$$
The Chern connection form rewritten in $x,s$ is

$$
\gamma = \frac{6 (b+2sb'+2+15sb^2+3s^2bb')d(s)}{ s(4+27sb^2)}-\frac{(2sb+4s^2b'+171sb^2+18s^2bb'+12b+24sb'+24)d(x)}{x(4+27sb^2)}.
$$
It is closed iff its coefficient by $\frac{d(x)}{x}$ is constant. Let us denote it by $L$ and rewrite this condition as an ODE for $b(s)$:
$$
(2sb'+b)(2s+9sb+12)=L(4+27sb^2).
$$
Substituting $\sqrt sb(s)=U(T), \ \ \sqrt s=T$ one arrives at the ODE for $U$
$$
[12+2T^2+9TU]\frac{dU}{dT}=L(4+27U^2).
$$
Solutions to this equation give the form 1), if $U$ is holomorphic at $T=0$, and the form 2), if $U$ has a pole at $T=0$.  The equation for $U$ is symmetric with respect to the involution $T\to -T$, $U\to -U$, hence one can choose $0 \leq \arg(U(0)) < \pi$. The detailed analysis is presented Appendix A.\\

{\it Case $a(s)=0$.} One easily checks that a suitable substitution (\ref{equivalence}) brings the equation to the form
$$p^3+y^4\beta (s)=0.$$ The Chern connection form of its web of solution   is
$$
\gamma = \frac{d\beta}{3\beta}+\frac{2}{3}\frac{\partial }{\partial_y}(\ln \beta )dy.
$$
Solving the equation $d(\gamma)=0$ for $\beta (s)$ one gets $\beta(s)=C_0s^k$, $C_0=const$. Note that $k$ should be integer nonnegative.
Applying a suitable substitution (\ref{equivalence}) we arrive at $p^3=y^4$. Thus, the general normal form is $p^3=x^{n}y^{4+m}$ with non-negative integer $n,m$. Applying Lemma \ref{hreduce} we prove that this form is equivalent to the form 2) with $n_0=2n+3$, $m_0=2m+1$, $L=-\frac{2}{3}$. Note, that it is sufficient to check the lemma hypothesis for the "basic" forms with  $m=n=0$.\\
\smallskip

{$\bullet$ \it Non-monic case.}
The function $F$ in (\ref{nonmonicCUB}) cannot vanish identically since 2 web directions cannot coincide identically. Thus its Tailor expansion has at least one term  and therefore the weights $w_x,w_y$ are rational. Introducing  a new invariant variable $s:=x^qy^{r-q}$, where $r,q$ are defined as in the monic case, we can rewrite (\ref{nonmonicCUB}) as
$$
\left(\frac{x}{y}\right)^2f(s)p^3+p+\left(\frac{y}{x}\right)g(s)=0,
$$
where $f(s) \not \equiv 0$. Then either $g(s)\equiv 0$ or the function $g$ can be brought to the form $g(s)=s^n$ by  some transformation (\ref{htransform}).

{\it Case $g(s)\equiv 0$.} By a suitable substitution (\ref{equivalence}) we can arrange $s=xy$. Then $d(\gamma)=0$ gives $f(s)=c_0s^c$. Applying now the inverse of substitution (\ref{equivalence}) we get the form 3).

{\it Case  $g(s)=s^n$.} By a suitable substitution (\ref{equivalence}) we can arrange $s=xy$ and $n=1$, which brings the equation into the form
$$
\left(\frac{x}{y}\right)^2f(s)p^3+p+y^2=0.
$$
In this equation the function $f$ must have the order at least 1 to give a non-singular coefficient by $p^3$ in the equation that we had before application of (\ref{equivalence}).
If $f(s)=sw(s)$, $w(0)\ne 0$ then substituting this expression into equation $d(\gamma)=0$ and analyzing its solution at $s=0$ with $w(0)\ne 0$ we get the following ODE for $w$:
$$
[2(1-s)+3s^2(3s-1)w]\frac{dw}{ds}=\frac{3}{2}[3s(3s+2)w+4]w,
$$
and the equation takes the form
$$
\left(\frac{2+m_0}{1+n_0}\right)^3x^{3+n_0}y^{m_0}w\left(\left[\frac{2+m_0}{1+n_0}\right]x^{1+n_0}y^{2+m_0}\right)p^3+p+x^{n_0}y^{3+m_0}=0.\\
$$
We claim that it is equivalent to 3) with the same values of parameters $n_0,m_0$. It is enough to apply Lemma \ref{hreduce} for "basic" equations with $n_0=m_0=0$.

Finally, let $f(s)=s^{2+l}v(s)$, where $l$ is non-negative integer and $v(0)\ne0$.
Substituting this expression into equation $d(\gamma)=0$ and analyzing its solution at $s=0$ with $v(0)\ne 0$ we get the following ODE for $v$:
$$
[2(1-s)+3s^{3+l}(3s-1)v]\frac{dv}{ds}=\frac{4+l}{2}[3s^{2+l}(2+3s)v+4]v.
$$
 with the corresponding normal form
 $$
 \left[\frac{1+m_0}{1+n_0}\right]^{4+l}x^{(1+n_0)(2+l)+2}y^{l+(l+2)m_0}{\textstyle v}\left(\left[\frac{1+m_0}{1+n_0}\right]x^{1+n_0}y^{1+m_0}\right)p^3+p+x^{n_0}y^{2+m_0}=0
 $$

We claim that for $l=0$ it is equivalent to 3), if  $n_3=2n_0+1$ and $m_3=2m_0$. Again it is sufficient to verify the hypothesis of Lemma \ref{hreduce} for the "basic" forms.
\hfill $\Box$\\

\subsection{Invariants}
Suppose equation (\ref{CUB}) is locally biholomorphic
to a weighted homogeneous one. How to determine the corresponding
normal form? There is a list of invariants that distinguishes between
the normal forms. Let $B_{q_0}$ be a small 4-dimensional ball over a singular point $q_0$ on the discriminant curve $\Delta$.
Obviously, the following objects are invariant under local biholomorphisms:
\begin{itemize}
   \item root multiplicity,
   \item projectivised  weights $[w_1:w_2]$,
   \item type of the discriminant curve singularity,
   \item periods of the  form $\gamma $ over cycles of the first homology group of $B_{q_0}\setminus\Delta$,
\end{itemize}
There is a subtler invariant. Consider the cross-ratio of the three web directions and the direction defined by the infinitesimal symmetry. This function is well-defined on  $B_{q_0} \setminus \Delta$ and is constant along the trajectories of the symmetry flow. Thus it is a function of the invariant parameter of the symmetry flow. For hyperbolic and elliptic weights this parameter is defined also at $q_0$. The limit value of this cross-ratio at $q_0$  is our fifth invariant. As the cross-ratio is dependent on the order of its arguments, we use the following  symmetrized form: multiply cubic form (\ref{binary}) with a 1-form vanishing on the trajectories of the symmetry group (for normal forms it could be $w_xxdy-w_yydx$), write the resulting quartic form $a_4dy^4+4a_3dy^3dx+6a_2dy^2dx^2+4a_1dydx^3+a_0dx^4$, compute $i:=a_0a_4-4a_1a_3+3a_2^2$ and  $j:=a_4a_2a_0+2a_1a_2a_3-a_2^3-a_4a_1^2-a_0a_3^2$. Then the invariant is $[i^3:j^2]$. The polynomials $i,j$ are well-known in the classical invariant theory, being invariants of the weights 4 and 6 respectively.

\begin{theorem}\label{invariants}
The normal forms described by Theorems \ref{parabolicTH},\ref{ellipticTH},\ref{hyperbolicTH} are pairwise not equivalent. The above defined five invariants distinguish between the equivalence classes.
\end{theorem}
{\it Proof:} In Appendix B, we present 3 tables describing invariants for parabolic, elliptic, and hyperbolic cases. The weights $w_x,w_y$ are defined (up to permutation) by the linear part of the symmetry operator thus giving also the type of the symmetry (parabolic, elliptic or hyperbolic). Then the root multiplicity $\mu$ eliminates the ambiguity in weight's order. Under the type of discriminant curve singularity we understand one of the following cases: non-singular variety, intersection (at $q_0$) of 2 or 3 non-singular varieties, intersection of 2 non-singular varieties and a singular variety passing through $q_0$. This information can be easily read off from the reduced equation of the discriminant curve at $q_0$. (Note that we do not need more subtle invariants like Tjurina number.)  The periods of the Chern connection form for the normal forms are determined by singular parts of $\gamma$ (i.e. by the equivalence class $[\gamma]$ of $\gamma$ modulo the subspace of holomorphic on $B_{q_0}$ forms) and by discriminant curve equations. The value of the cross-ratio invariant is presented for indicated parameters to distinguish between the forms when other invariants are not effective. More details are given in Appendix B. \hfill $\Box$\\

\section{Hexagonal 3-webs with holomorphic Chern connection.}
The obtained classification of weighted homogeneous ODEs with a hexagonal web of solutions allows to classify 3-webs with an exact Chern connection and an infinitesimal symmetry vanishing at the singular point.
\begin{proposition}\label{exactvanishing euler}
If an infinitesimal symmetry of equation (\ref{binary}) vanishes at $(0,0)$ and the Chern connection is exact, then the equation  is equivalent
to a weighted homogeneous one and the symmetry to an Euler vector field.
\end{proposition}
{\it Proof:}
In fact, choosing the first integrals as in Lemma \ref{first integrals for exact} we have $c_1=c_2=c_3=0$ in formula (\ref{Xu}). For example, $\left .X(u_1)\right |_{0}=C\left .u_1\right |_0+c_1=c_1$. On the other hand $X(u_1)=k(p_2-p_3)(\eta -p_1\xi)=0$. Whence $c_1=0$. Therefore $C\ne 0$ and the equation is weighted homogeneous due to Theorem \ref{Xeuler}. \hfill $\Box$\\
\begin{corollary}
If an infinitesimal symmetry of equation (\ref{binary}) vanishes at the singular point and the Chern connection is exact, then the squares of abelian integrals  $u_i$ satisfy equation (\ref{cub first integrals}).
\end{corollary}
The following classification is an immediate consequence of the above proposition and Theorems \ref{parabolicTH},\ref{ellipticTH},\ref{hyperbolicTH},\ref{invariants}.
\begin{theorem}\label{homoexact}
Suppose ODE (\ref{binary}) admits an infinitesimal symmetry $X$ vanishing at the point $(0,0)$ on the discriminant curve $\Delta$ and the germ of the Chern connection form is exact  $\gamma=d(f)$,
where  $f$ is some function germ. Then  the equation and the symmetry are equivalent to one of the
following normal forms:
$$
\begin{array}{lll}
1) & y^{m_0}p^3-p=0, & X=(2+m_0)x\partial _x+2y\partial _y,\\
2) & p^3+2xp+y=0,    & X=2x\partial _x+3y\partial _y,\\
3) & (p-\frac{2}{3}x)(p^2+\frac{2}{3}xp+y-\frac{2}{9}x^2)=0, & X=x\partial _x+2y\partial _y,\\
4) & p^3+4x(y-\frac{4}{9}x^3)p+y^2+\frac{64}{81}x^6-\frac{32}{9}yx^3=0, & X=x\partial _x+3y\partial _y,\\
5)& p^3+xy^2p+\frac{2}{\sqrt{27}}\frac{x^{\frac{3}{2}}y^3}{\tan(\frac{4}{\sqrt{3}}x^{\frac{3}{2}})}=0, & X=y\partial _y, \\
6) & p^3+y^2p=\frac{2}{\sqrt{27}}y^3\tan(2\sqrt{3}x+L), & X=y\partial _y, \\
7)&p^3+y^{3+m_0}p+y^{\frac{9+3m_0}{2}}U\left(\left[(m_0+1)\right]xy^{\frac{1+m_0}{2}}\right)=0, & X=(1+m_0)x\partial _x-2y\partial _y,\\
8)&p^3+xy^{3+m_0}p-\frac{x^{\frac{3}{2}}y^{\frac{9+3m_0}{2}}}{V\left(\left[\frac{2}{3}(m_0+1)\right]x^{\frac{3}{2}}y^{\frac{1+m_0}{2}}\right)}=0, & X=(1+m_0)x\partial _x-3y\partial _y,\\
\end{array}
$$

\noindent The function $U(T)$  is defined, with $L=-\frac{2(m_0+3)}{(m_0+1)}$ and suitable constants $C_1,C_2$, by the relations
$$
\begin{array}{ll}
\frac{T}{3\sqrt3L}=\frac{f'\left(-\arctan\left(\frac{3\sqrt3}{2}U\right)\right)}{f\left(-\arctan\left(\frac{3\sqrt3}{2}U\right)\right)},& f(z)=\cos^{-\mu} (z)[C_1P^{\mu}_{\nu}(\sin z)+ C_2Q^{\mu}_{\nu}(\sin z)]\\
f\left(-\arctan\left(\frac{3\sqrt3}{2}U(0)\right)\right)\ne 0, & f'\left(-\arctan\left(\frac{3\sqrt3}{2}U(0)\right)\right)=0.\\
\end{array}
$$
The initial value of $U$ vanishes $U(0)=0$ if at least one of the numbers $n_0$, $m_0+1$ is odd. If $U(0)\ne 0$ one can choose $0 \leq \arg(U(0)) < \pi$.

The function  $V(T)$ is defined, with $L=-\frac{5m_0+17}{3(m_0+1)}$, by the relations
$$
\begin{array}{lll}
\frac{1}{3\sqrt3 L}T=\frac{f'\left(-\arctan\left(\frac{2}{3\sqrt3}V\right)\right)}{f\left(-\arctan\left(\frac{2}{3\sqrt3}V\right)\right)},& f(z)=\sin^{\mu} (z)P^{\mu}_{\nu}(\cos z).& \\
\end{array}
$$
In the above formulas, $P^{\mu}_{\nu}(z),Q^{\mu}_{\nu}(z)$ are Legendre`s functions for $\mu=\frac{1}{2}(1-\frac{1}{3L})$, $\nu=\frac{1}{2}(\frac{1}{L}-1)$,  $m_0$ is non-negative integer, for the form 6) with $L\ne 0$ one can choose  $0 \leq \arg(L) < \pi$.\\
The weights $[w_1:w_2]$, the root multiplicity and the invariant $[i^3:j^2]$ uniquely determine  the normal form.
\end{theorem}
To complete classification of all implicit ODEs with exact Chern connection and with at least one infinitesimal symmetry, we have to consider the case when the symmetry does not vanish.
Along with Lemma \ref{productLM} we will need the following one.
\begin{lemma}\label{ratioLM}
Suppose $f:(\mathbb C,0)\to (\mathbb C,f_0)$ is holomorphic, $f_0:=f(0)\ne 0$, $f_{n-1}:=\left.\frac{1}{n!}\frac{df^{n-1}}{dx^{n-1}}\right|_{t=0}$. Then the following ODEs have analytic solutions\\ $\psi: (\mathbb C,0) \to (\mathbb C,0)$ with $\left.\frac{d\psi}{dt}\right|_{t=0}\ne 0$
\begin{enumerate}
\item  $\frac{f(\psi)}{\psi }\frac{d\psi}{dt}=\frac{f_0}{t}$,
\item  $\frac{f(\psi)}{\psi ^q}\frac{d\psi}{dt}=\frac{f_0}{t^q}$, where $q\in \mathbb Q,\ q\notin \mathbb N,\ q>0$,
\item  $\frac{f(\psi)}{\psi ^n}\frac{d\psi}{dt}=\frac{f_0}{t^n}+\frac{f_{n-1}}{t}$, where $n\in \mathbb N,\ n>1$.
\end{enumerate}
\end{lemma}
{\it Proof:} Expanding $f$ in Taylor series $f(\psi)=\sum_{k=1}^{\infty} f_k\psi ^k$ and integrating we get  $\ln \psi +\tilde{f}(\psi)=\ln t$ for the first equation, where the function $\tilde{f}$ is analytic. Therefore
$\psi e^{\tilde{f}(\psi)}=t$ gives the desired solution. Integrating the second equation we arrive at $\psi \tilde{f}(\psi)=t$, where the function $\tilde{f}$ is analytic and $\tilde{f}(0)\ne 0$. In the third case we have the following relation after integration
$$
\frac{\frac{f_0}{1-n}\tilde{f}(\psi)+f_{n-1}\psi ^{n-1}\ln \psi}{\psi^{n-1}}=\frac{\frac{f_0}{1-n}+f_{n-1}t^{n-1}\ln t}{t^{n-1}},
$$
where the function $\tilde{f}$ is analytic and $\tilde{f}(0)=1$. Substituting $\psi=t(1+z)$ we see that the terms with $\ln t$ are canceled and the resulting equation could be locally resolved for $z$ at the point $(t,z)=(0,0)$.
\hfill $\Box$\\

\begin{theorem}\label{shiftTH}
 Suppose  ODE (\ref{binary}) admits an infinitesimal symmetry that does not vanish at the point $(0,0)$ on the discriminant curve $\Delta$ and the germ of the Chern connection form is exact  $\gamma=d(f)$,
where  $f$ is some function germ. Then for some natural $m_0$  the equation is equivalent to the
following one
$$
 y^{m_0}p^3-p=0,
$$
admitting the 2-dimensional symmetry algebra generated by the operators\\ $\{\partial _x,\ (2+m_0)x\partial _x+2y\partial _y\}$.
\end{theorem}
{\it Proof:} If the infinitesimal symmetry does not vanish at $(0,0)$ it defines a direction $\tau$. The relative position of the direction $\tau$ with respect to the web directions $p_1,p_2,p_3$ could be one of the following:
\begin{enumerate}
\item $\tau$ does not coincide with any $p_i$,
\item $\tau$ coincides with the double direction $p_1=p_2\ne p_3$,
\item  $\tau$ coincides with the simple direction $p_3$, where $p_3\ne p_1=p_2$,
\item $\tau$ coincides with the triple direction $p_1=p_2=p_3$.
\end{enumerate}
{$\bullet$ \it Case 1.} First let us bring the infinitesimal symmetry to the form $X=\partial _y$ by a suitable coordinate transform. Then $dx=0$ is not a root of equation (\ref{binary}) hence $K_3(0,0)\ne 0$. Thus the corresponding cubic equation is monic, i.e. it has the form (\ref{general_cub}). Due to the symmetry, the functions $a,b,c$ do not depend on $y$. The biholomorphisms respecting $\partial _y$ have the form
\begin{equation}\label{preserveDY}
y=c\tilde{y}+\varphi(\tilde{x}), \ \ x=\psi(\tilde{x}).
\end{equation}
Choosing $\varphi$ to satisfy $3\varphi _x(x)+a(x)=0$ we kill the coefficient by $p^2$. Now our equation has the form (\ref{CUB}). If $A(x)\equiv 0$ and $B(x)=x^n\beta(x)$ with $\beta(0)\ne 0$, we use a biholomorphism of the form (\ref{preserveDY}) with $\varphi=0$ (and Lemma \ref{productLM}) to bring our equation to the normal form $p^3+x^{n}=0$, whose Chern connection $\gamma =ndx/3x$ is closed but not exact for $n > 0$.
If $A(x)=x^n\alpha(x)$ with $\alpha (0)\ne 0$ we again use a biholomorphism of the form (\ref{preserveDY}) with $\varphi=0$  to bring our equation to the form $p^3+x^np+B(x)=0$.
The substitution (\ref{equivalence}) linear in $y$ reduces it to $p^3+p+b(x)=0$ with the Chern connection $$\gamma =\frac{9bb_xdx+6b_xdy}{4+27b^2}.$$ Now hexagonality is equivalent to the following ODE for $b$:
$$\frac{6b_x}{4+27b^2}=const.$$ Its general solution is $b(x) = \frac{2}{3\sqrt3}\tan(Lx+L_1)$. The corresponding Chern connection
$$
\gamma=\frac{L}{3}\{\tan(Lx+L_1)dx+\sqrt3dy\}
$$
is exact but the   discriminant  $D=-4(1+\tan(Lx+L_1)^2)$  never vanishes. Applying an inverse substitution $x\to x^{1+\frac{n}{2}}$ we get a non-exact Chern connection with the term $ndx/2x$. Thus the case 1 does not give singular webs with exact Chern connection.\\
{$\bullet$ \it Case 2.} We adjust coordinates so that:
\begin{enumerate}
\item the simple web direction coincides with that of the $x$-exes,
 \item the double direction is that of the $y$-axes,
  \item$X=\partial _y$.
\end{enumerate}
   Since one of the web directions is simple and transverse to $X$, the corresponding foliation is defined by an analytic ODE $\frac{dy}{dx}=p(x)$, where $p$ does not depend on $y$ due to the symmetry. Choosing $\tilde{y}=y-\int p(x)dx$ we bring this foliation to the form $\tilde{y}=const$ and preserve $X$. Thus we can assume that one of the roots of our cubic ODE vanishes identically, i.e., the equation is $p\left(K(x)p^2+L(x)p+M(x)\right)=0$ with $M(0)\ne 0$ as the root $p=0$ is simple. Therefore our equation can be written as $p(K(x)p^2+L(x)p+1)=0$. Further, the direction $dx=0$ has multiplicity two, hence $K(0)=L(0)=0$. Using a biholomorphism of the form (\ref{preserveDY})  with $\varphi=0$ (and Lemma \ref{ratioLM}) we reduce the equation to the form
\begin{equation}\label{Kdelta}
K(x)p^3+\frac{x^{n+1}}{1+\delta x^{n}}p^2+p=0,
\end{equation}
where $n$ is a nonnegative integer and $\delta =0$ or $\delta=1$. (One can always make $f_0=1$ and $f_{n-1}=1$ or $f_{n-1}=0$ by a suitable scaling of $x$ and $y$.)

If $n=0$ then the equation is $K(x)p^3+xp^2+p=0$. Now the condition $d(\gamma)=0$ gives $K(x)=\frac{x^2}{4}(1+4L_1x^{4L})$. The corresponding Chern connection $\gamma =(1+2L)\frac{dx}{x}+Ldy$ is exact iff $L=-\frac{1}{2}$. Then $L_1=0$ as  $dx=0$ is one of the web directions (namely, of multiplicity two). Therefore our equation $p(xp+2)^2=0$ has two identically coinciding roots.

If $n>1$ then the substitution (\ref{equivalence}) linear in $y$ reduces the equation to $K(x)p^3+\frac{x^2}{1+\delta x}p^2+p=0$.
For $\delta =0$ we have $K(x)=\frac{x^4}{4}(1+L_1e^{\frac{4L}{x}})$. For equation (\ref{Kdelta}) to be holomorphic at $(0,0)$ it is necessary that either $L_1=0$ or $L=0$. The condition $L_1=0$ implies that the equation has two identically coinciding roots. If $L=0$ then the Chern connection is not exact having the term $dx/x$.\\
{$\bullet$ \it Case 3.} We adjust coordinates so that:
\begin{enumerate}
\item the simple web direction coincides with that of the $x$-exes,
 \item the double direction is that of the $y$-axes,
  \item$X=\partial _x$.
\end{enumerate}
The biholomorphisms respecting $\partial _x$ have the form
\begin{equation}\label{preserveDX}
x=c\tilde{x}+\varphi(\tilde{y}), \ \ y=\psi(\tilde{y}).
\end{equation}
As $dy=0$ is a simple web direction for equation (\ref{binary}) we have $K_0(0)=0$, $K_1(0)\ne 0$. (Observe that all $K_i$ depend only on $y$.)  Thus we can divide the equation by $K_1$ and assume $K_1\equiv 1$. Further, since $dx=0$ is the double web direction we have $K_3(0)=K_2(0)= 0$. Choosing $\psi(y)=y$ and $\varphi$ satisfying $K_2(y)+2\varphi _y(y)+3K_0(y)\varphi _y(y)=0$  we kill the coefficient by $p^2$. Therefore our equation can be written as $K(y)p^3+p+N(y))=0$.

If $N\equiv 0$ we use a biholomorphism of the form (\ref{preserveDX})  with $\varphi\equiv 0$ and Lemma \ref{productLM} to reduce the equation to the form $y^{n}p^3-p=0$, which is already known to have the desired properties. (See Theorem \ref{homoexact}.)

If $N(y)=y^n\beta(y)$, $\beta(0)\ne 0$ we use a biholomorphism of the form (\ref{preserveDX})  with $\varphi\equiv 0$ and Lemma \ref{ratioLM} to reduce the equation to the form
\begin{equation}\label{Kyn}
K(y)p^3+p+y^n,
\end{equation}
where $n$ is natural.
If $n>1$ then the substitution (\ref{equivalence}) linear in $x$ reduces the equation to $K(y)p^3+p+y^2$.
The equation $d(\gamma)=0$ implies
$$
K(y)=\frac{4}{y^4\left(L_1 e^{-2 \frac{L}{y}}-27\right)}.
$$ Therefore equation (\ref{Kyn}) can not be holomorphic at $(0,0)$ for any choice of $L_1$ and $L$.

If $n=1$ then the equation $d(\gamma)=0$ implies
$$
K(y)=\frac{4}{-27y^2+4L_1y^{2+2L}}.
$$
This function is holomorphic and vanishing at $0$ iff $2+2L=-n$, where $n$ is natural. Thus our equation is
$$
\frac{4y^n}{4L_1-27y^{2+n}}p^3+p+y, \ \  \ L_1\ne 0.$$
We claim that this equation is locally equivalent to   $y^{n}p^3-p=0$. To prove this we show that except for the symmetry $X=\partial _x$ the equation admits a symmetry vanishing at $(0,0)$ and therefore is equivalent to one from the list of Theorem \ref{homoexact}. The point $(0,0)$ is a singular point of the discriminant curves of the second and the third form hence  they can not admit a non-vanishing symmetry. The last two forms of the list have a triple root at $(0,0)$. Thus the equation is equivalent to the first form.

The components $\xi,\eta$ of the infinitesimal symmetry $X$ satisfy so-called defining Lie equations (see \cite{Lds}). To write them for our particular case
one has to extend the infinitesimal symmetry action on $p$ by the formula
$$
\tilde{X}=X+\zeta(x,y,p)\partial _p,
\ \  {\rm where} \ \ \zeta(x,y,p)=D\eta (x,y)-pD\xi(x,y), \
\ D=\partial _x +p\partial _y,
$$
(the extended action respects the contact field $dy-pdx=0$ at $\mathbb C ^2 \times \mathbb P^1(\mathbb C)$),
apply the symmetry operator to the equation and get a polynomial of second degree in $p$ as the rest by division by the cubic equation. Splitting this polynomial with respect to $p$ we obtain  three equations. From these equations we have:
$$
\eta _y= \xi_x-\frac{4nL_1+(36-9n)y^{2+n}}{2y(4L_1-27y^{2+n})}\eta, \ \ \xi _y = \frac{3(2+n)}{4y^{-n}L_1-27y^2}\eta, \ \ \eta _x=-\frac{1}{2}(2+n)\eta .
$$
From the last equation it is immediate that $\eta= e^{-(1+\frac{n}{2})x}F(y)$. Substituting this expression into the first of the above equations one concludes that
$\xi _x= e^{-(1+\frac{n}{2})x}G(y)$. Hence $\xi = \frac{e^{-(1+n/2)x}}{-1-n/2}G(y)+C$. (The constant $C$ here corresponds to the symmetry $\partial _x$.) Moreover, compatibility conditions imply the following equations for the functions $F,G$:
$$
G_y= -\frac{3(2+n)^2F}{2(4y^{-n}L_1-27y^2)},\ \ \ F_y= \frac{(-4nL_1-36y^{2+n}+9ny^{2+n})F+(8yL_1-54y^{3+n})G}{2y(4L_1-27y^{2+n})}.
$$
Solving the first equation for  $F$:
\begin{equation}\label{expresF}
F= -\frac{2(4L_1-27y^{2+n})}{3(2+n)^2}\frac{G_y}{y^{n}}
\end{equation}
and substituting this expression into the second we get a second order ODE for $G$. Fortunately, it can be integrated in closed form:
$$
G(y)=C_1\left(\frac{2\sqrt L_1}{3\sqrt 3}-y^{1+\frac{n}{2}}\right)^{\frac{1}{3}}+C_2\left(\frac{2\sqrt L_1}{3\sqrt 3}+y^{1+\frac{n}{2}}\right)^{\frac{1}{3}}
$$
Choosing $C_1=C_2=1$ we get an even analytic function of $y^{1+\frac{n}{2}}$, i.e., analytic in $y$. Moreover, the first nonconstant term in the Taylor expansion of $G$ is $y^{2+n}$. Therefore the function $F$ defined by (\ref{expresF}) is holomorphic at $0$ and $F(0)=0$. Selecting the constant $C$ in the formula for $\xi$ we obtain a symmetry vanishing at $0$.
\smallskip

{$\bullet$ \it Case 4.} We adjust coordinates so that:
\begin{enumerate}
\item the triple web direction coincides with that of the $x$-exes,
  \item$X=\partial _x$.
\end{enumerate}
We choose the abelian first integrals $u_i$ to vanish at $(0,0)$, which implies $c_i=0$ in formula (\ref{Xu}).
Therefore $u_i(x,y)=k_ie^{C\left(x-\int \frac{dy}{p_i(y)}\right)}$ for some constants $k_i$.
In adapted coordinates the equation has the form
$$
p^3+a(y)p^2+b(y)p+c(y)=0,
$$
where $a(0)=b(0)=c(0)=0$.
Its roots can be expanded in Puiseux series $p_i(y)=y^{q_i}\rho_i(y^{\frac{1}{3}})$, where $q_i>0$ are rational and $\rho _i,$ with $\rho _i(0)\ne 0$ are analytic.

Due to Lemma \ref{first integrals for exact} the functions $u_i$ are algebraic over $\mathcal{O}_0$, consequently we have $q_i\le 1$. If $q_i<1$ then $u_i$ does not vanish at $(0,0)$. Whence $q_i=1$. Using Puiseux expansions of  $p_i(y)$ we obtain  $u_i(x,y)=k_iy^{r^i}e^{C\left(x-f_i(y^{1/3})\right)}$, where the functions $f_i$ are analytic and satisfy $f_i(0)=0$. The condition $$u_1+u_2+u_3\equiv 0$$ implies 1) $r_1=r_2=r_3$ and 2) $f_1(t)\equiv f_2(t)\equiv f_3(t)$. But that means that all 3 roots coincide identically. Thus in this case we also do not obtain new forms. \hfill $\Box$\\
\smallskip

Note that if the symmetry algebra is 3-dimensional and the Chern connection form is exact, then the roots of equation (\ref{binary}) are simple. In fact, there are two symmetries $X_1$, $X_2$ satisfying $X_i(u_j)=\delta_{ij}$. The function $k$ in equations (\ref{du}) for abelian first integrals can be reduced to $k=1$ (Lemma \ref{k1}).
\begin{corollary}
Suppose equation (\ref{binary}) has a non-trivial symmetry algebra at $(0,0)$ and the Chern connection of the web germ of solutions is exact.
\begin{itemize}
\item If the symmetry algebra is 3-dimensional, then the equation has  simple roots at $(0,0)$ and the web germ is regular.
 \item If the symmetry algebra is 2-dimensional, then the equation has a double root at $(0,0)$ and is equivalent to the form  1) in Theorem \ref{homoexact}.
 \item If the symmetry algebra is 1-dimensional, then the equation is equivalent to one of the other forms in Theorem \ref{homoexact}.
\end{itemize}
\end{corollary}

\noindent {\bf Remark.} Observe that in the proof of Theorem \ref{shiftTH} the first 3 cases can be considered in full generality to obtain normal forms without the condition of analyticity of the closed Chern connection (at least in terms of fixed solution of some ODE, like the form 6) in Theorem \ref{homoexact}). In the last case the difficulties are much stronger: one has to consider singular systems of ODEs instead of one singular ODE in Lemmata \ref{productLM} and \ref{ratioLM}.

\section{Concluding remarks}
\subsection{Geometric construction for characteristic webs}
There is a geometric construction for characteristic webs on solutions of associativity equations (\ref{ass}) and (\ref{ass2}) (see \cite{Aw}).  Suppose the solution describes a semi-simple Frobenius manifold. Then at each regular point the tangent space is decomposed in the direct sum of 3 one-dimensional algebras $\mathbb C$. Take a vector field $e_i$ corresponding to one of the three idempotents (the unities of these one-dimensional algebras) and the unity vector field $e$. These 2 vector fields $e,e_i$ define a two-dimensional integrable distribution. Thus we have 3 foliations, one for each $i=1,2,3.$ Choose a surface transverse to the field $e$.    Then the leaves of these 3 foliations cut a 3-web on the surface. This web is hexagonal; it follows from existence of {\it local canonical coordinates} (see \cite{Dga},\cite{Mf}), which are sometimes called also {\it Dubrovin coordinates}.  This 3-web is equivalent to the characteristic web.

\subsection{Frobenius 3-folds from 3-webs} The natural question is, in what extent the above construction can be inverted to recover a Frobenius  3-fold germ  starting with a singular web with infinitesimal symmetry and holomorphic Chern connection? There are good chances. We have a symmetry, which is equivalent to an Euler vector field. Therefore we have also good candidates for flat coordinates. Finally, the web directions suggest idempotent directions. The details are discussed in \cite{Af}, where the associativity equations were interpreted geometrically in terms of the web theory. There is a hope that a similar interpretation is possible for any dimension.

\subsection{Generalization}
The above geometric construction can be generalized to higher dimensions. As a result we obtain, for $n$-dimensional Frobenius manifold, a collection of $n$ commuting vector fields $v_i$ in $(\mathbb C^{n-1},0)$, satisfying the equation $\sum^n_{i=1}v_i=0$, i.e., a "flat" $n$-web germ of curves in $(\mathbb C^{n-1},0)$ admitting a "linear" symmetry. We hope that using the language of the web theory will bring better insight into the structure of the discriminant set of Frobenius manifolds.

\subsection{Acknowledgement}
The author thanks the hospitality of the Institute of Mathematical and Computer Sciences of S\~ao Paulo University USP-ICMC in S\~ao Carlos, where this study was initiated, and M.A.S.Ruas in particular. This
research was partially supported by CNPq grant 454618/2009-3 and by the National Institute of Science
and Technology of Mathematics INCT-Mat.

\section{Appendix A: solutions to $[12+2t^2+9tU]\frac{dU}{dt}=L(4+27U^2)$}
First consider the case of a holomorphic solution $U$.  The substitution
$$
T=3\sqrt 3 L\frac{f'(z)}{f(z)}, \ \ \ U=-\frac{2 }{3\sqrt 3}\tan(z),
$$
linearizes the problem; the function $f$ satisfies the linear ODE
$$
f''-\frac{\tan(z)}{3L}f'+\frac{2}{9L^2}f=0,
$$
whose general solution is expressed in terms of the Legendre functions.
In fact, the substitution $f(z)=(1-x^2)^{-\frac{\mu}{2}}g(x)$, $x=\sin z$ transforms the above equation to the Legendre equation
\begin{equation}\label{Legendre}
(1-x^2)g''-2xg'+\left[\nu(\nu+1)-\frac{\mu^2}{1-x^2} \right]g=0
\end{equation}
with $\mu=\frac{1}{2}(1-\frac{1}{3L})$, $\nu=\frac{1}{2}(\frac{1}{L}-1)$. Thus we obtain the normal form 1) with $U$ defined by (\ref{UPQ}).

It is easily seen that that $U(T)$ is allowed to have the pole of order 1, which corresponds to $n_0\ge 1$ in the form with $U$, or of order 2, which corresponds to $n_0\ge 4$. The substitution
$$
T=3\sqrt 3 L\frac{f'(z)}{f(z)}, \ \ \ U=\frac{2}{3\sqrt 3 \tan(z)},
$$
now brings the equation for $U$ to
\begin{equation}\label{f_for_V}
f''+\frac{1}{3L\tan(z)}f'+\frac{2}{9L^2}f=0.
\end{equation}
The singular point $z=0$ is regular. A standard local analysis of solutions (see, for example, \cite{Gk}) gives the following types of solutions with analytic non-vanishing at $z\ne 0$ functions $A,B,C$ and a constant $\lambda$.

1. If $\rho:=1-\frac{1}{3L}\ne \mathbb Z$ then $f(z)=c_1A(z)+c_2z^{\rho}B(z)$.

2. If $\rho =-n,$ $n\in \mathbb N$ then $f(z)=c_1A(z)+c_2[z^{-n}B(z)+\lambda \ln z A(z)]$.

3. If $\rho =0$  then $f(z)=c_1A(z)+c_2[\kappa \ln z A(z)+z \psi(z)]$, where $\kappa \ne 0$ and  $\psi$ is analytic.

4. If $\rho =n,$ $n\in \mathbb N$ (note that $n>1$) then $f(z)=c_1z^nA(z)+c_2[C(z)+\lambda z^n \ln z B(z)]$.

\noindent The function $U$ has pole of order 1 at $z=0$ iff $f(z)$ is analytic with $f(0)\ne 0$, $f'(0)=0$, $f''(0)\ne 0$. If there is an analytic non-vanishing at $z=0$ solution, then it automatically verifies $f'(0)=0$, $f''(0)\ne 0$. A solution of types 1,2 or 3 suits iff $c_2=0$, thus giving (\ref{VPQ}), where we use $V=-\frac{1}{U}$. In fact,
the substitution $f(z)=(1-x^2)^{\frac{\mu}{2}}g(x)$, $x=\cos z$ transforms equation (\ref{f_for_V}) to the Legendre equation (\ref{Legendre})
with the same $\nu$ and $\mu$. For the solutions of the type 4, an analysis of series expansions of the functions $P^{\mu}_{\nu}(z),Q^{\mu}_{\nu}(z)$ at $z=1$ (see, for instance, \cite{Eh}), shows that  $\alpha =0$ always holds for odd $n$ and $\alpha \ne 0$ always holds for even $n$. Therefore a solution with desired properties exists only for odd $n$ and is of the form
\begin{equation}\label{fullPQ}
f(z)=\sin^{\mu} (z)[C_1P^{\mu}_{\nu}(\cos z)+ C_2Q^{\mu}_{\nu}(\cos z)]
\end{equation}
where $C_1,C_2$ are chosen to guarantee $c_2\ne 0$. Lemma \ref{hreduce} implies that all corresponding ODEs are equivalent, thus one can choose  $C_2=0$ and get (\ref{VPQ}).

\noindent The function $U$ has pole of order 2 at $z=0$ iff $f(z)$ is of the type 1 with $c_1\ne 0$, $c_2\ne 0$. Therefore $L= -\frac{2}{3}$ and
the solution is of the form (\ref{fullPQ}).  Due to Lemma \ref{hreduce} all corresponding ODEs are equivalent and we can choose  $C_1=C_2=1.$ This gives (\ref{VPQ2}).

\section{Appendix B: tables of invariants}

\noindent Here is the table for the parabolic case.
\smallskip

\hspace{-0.4cm}\begin{tabular}{|c|c|c|c|c|l|}
  \hline
\# &$\mu$&$[w_1:w_2]$&$\Delta$      &$[\gamma]$ & $[i^3:j^2]$ \\
  \hline
1 &  3   &$[0:1]$     &$yx^N=0$ &$\frac{N}{3}\frac{dx}{x} + 2\frac{dy}{y}$     &   $[0:1]$   \\
  \hline
2 &  3   &$[0:1]$     &$yx^N=0$     &$\frac{N}{2}\frac{dx}{x}+(2-\frac{L_0(N+2)}{2\sqrt{3}})\frac{dy}{y}  $     &    $[1:\frac{\tan^2(L_1)}{-27}]$    \\
  \hline
3 &  3   &$[0:1]$     &$yx^N=0$       &$\frac{N}{3}\frac{dx}{x} + (2-\frac{L_0(N+3)}{2\sqrt{3}})\frac{dy}{y}$& $[0:1]$ \\
   \hline
4 &  2   &$[0:1]$    &$x^2y^{N}=0$&$\frac{dx}{x}$      &  \\
     \hline
\end{tabular}\\

\smallskip

\noindent The table for
the  hyperbolic case is as follows:

\smallskip

\hspace{-0.4cm}\begin{tabular}{|c|c|c|c|c|l|}
  \hline
\# &$\mu$&$[w_1:w_2]$&$\Delta$      &$[\gamma]$ & $[i^3:j^2]$ \\
  \hline
1 &  3   &$[-(m_0+1):n_0+2]$     &$yx^{n_0}=0$        &$ \frac{n_0}{2}\frac{dx}{x}+\{2+(\frac{L}{2}+1)(1+m_0)\}\frac{dy}{y}  $     &    $[1:\frac{U^2(0)}{-4}]$    \\
  \hline
2 &  3   &$[-(m_0+1):n_0+3]$     &$yx^{n_0}=0$        &$ \frac{n_0}{3}\frac{dx}{x}+\{2+(\frac{L}{2}+\frac{5}{6})(1+m_0)\}\frac{dy}{y}  $     &    $[0:1]$    \\
  \hline
2.2  &  3   &$[-(m_0+1):n_0+3]$     &$yx^{n_0-3}=0$        &$ \frac{n_0-3}{6}\frac{dx}{x}+\frac{m_0+7}{3}\frac{dy}{y}  $     &    $[0:1]$    \\
  \hline
3 &  2   &$[-(m_0+2):n_0+1]$     &$xy^{m_0}=0$       &$\frac{n_0+3}{2}\frac{dx}{x}$                                                & \\
  \hline
4 &  2   &$[-(m_0+1):n_0+1]$    &$xy=0$&$(1+\frac{(n_0+1)(l_0+3)}{2})\frac{dx}{x}+(1+l_0)(1+m_0)\frac{dy}{y}$      &  \\
     \hline
\end{tabular}\\

\smallskip

\noindent Some comments on the hyperbolic case:
to $U\equiv0$ corresponds $L=0$, the case 2.2 corresponds to $L=-\frac{2}{3}$.\\
\noindent In the table for elliptic case we add the exponents $\alpha,\beta$ generating normal forms from the corresponding "basic" equation of the list of Theorem \ref{ellipticTH}.\\
Some comments:\\
  1) If $n_0$ or $m_0$ comes with the negative sign in formulas for $\alpha$ or $\beta$ then it is either $0$ or $1$.\\
  2)  $\lambda$ is a non-vanishing constant (its value can be easily computed).\\
  3)  The invariant $[i^3:j^2]$ is used only once to distinguish between the forms 18) and 26).

\bigskip

\hspace{-1.0cm}\begin{tabular}{|c|c|c|c|c|c|}

  \hline
\#&$\alpha,\beta$                     &$\mu$&$[w_1:w_2]$    &$\Delta$        &$[\gamma]$     \\
  \hline
1 &$\scriptstyle 1-\frac{n_0}{2}, 1+\frac{m_0}{2}$ &  2  &$\scriptstyle[m_0+2:2-n_0]$&$\scriptstyle x^{n_0}y^{m_0}$&$\frac{n_0}{2}\frac{dx}{x}$  \\
  \hline
1 &$\scriptstyle1+\frac{n_0}{2}, 1-\frac{m_0}{2}$ &  3  &$\scriptstyle[2-m_0:n_0+2]$&$\scriptstyle x^{n_0}y^{m_0}$&$\frac{n_0}{2}\frac{dx}{x}+\frac{m_0dy}{y}$\\
  \hline
2 &$\scriptstyle1+\frac{n_0}{3}, 1-\frac{m_0}{2}$ &  3  &$\scriptstyle[2-m_0:n_0+3]$&$\scriptstyle x^{n_0}y^{m_0}(4x^{3+n_0}+27y^{2-m_0})$&$\frac{n_0}{3}\frac{dx}{x}+\frac{m_0dy}{y}+\frac{d\ln (4x^{3+n_0}+27y^{2-m_0})}{2}$ \\
   \hline
3 &$\scriptstyle1+\frac{n_0}{3}, 1-\frac{m_0}{2}$ &  3  &$\scriptstyle[2-m_0:n_0+3]$&$\scriptstyle x^{n_0}y^{m_0}(32x^{3+n_0}+27y^{2-m_0})$&$\frac{n_0}{3}\frac{dx}{x}+\frac{m_0dy}{y}$  \\
  \hline
4 &$\scriptstyle1+\frac{n_0}{2}, 1$               &  3  &$\scriptstyle[1:n_0+2]$    &$\scriptstyle yx^{n_0}(12y+\lambda x^{2+n_0})$          &$\frac{n_0}{2}\frac{dx}{x}+\frac{2}{3}\frac{dy}{y}+\frac{d\ln (12y+\lambda x^{2+n_0})}{2}$     \\
\hline
5 &$\scriptstyle1+\frac{n_0}{2}, 1$               &  3  &$\scriptstyle[1:n_0+2]$    &$\scriptstyle yx^{n_0}(3y+\lambda x^{2+n_0})$          &$\frac{n_0}{2}\frac{dx}{x}+\frac{2}{3}\frac{dy}{y}$   \\
\hline
6 &$\scriptstyle1+\frac{n_0}{2}, 1$               &  3  &$\scriptstyle[1:n_0+2]$    &$\scriptstyle x(x^{2+n_0}-6\lambda y)(x^{2+n_0}-3\lambda y)$          &$\frac{2n_0+1}{3}\frac{dx}{x}+\frac{d\ln (x^{2+n_0}-6\lambda y)}{2}$      \\
\hline
7 &$\scriptstyle1+\frac{n_0}{2}, 1$               &  3  &$\scriptstyle[1:n_0+2]$    &$\scriptstyle x^{n_0}(x^{2+n_0} -12\lambda y)(x^{2+n_0}-3\lambda y)$  &$\frac{n_0}{2}\frac{dx}{x}+\frac{d\ln (x^{2+n_0}-12\lambda y)}{2}$        \\
     \hline
8 &$\scriptstyle1+\frac{n_0}{2}, 1$               &  3  &$\scriptstyle[1:n_0+2]$    &$\scriptstyle x^{n_0}(x^{2+n_0} -3\lambda y)(2x^{2+n_0}+3\lambda y)$  &$\frac{n_0}{2}\frac{dx}{x}$    \\
     \hline
9 &$\scriptstyle1+\frac{n_0}{2}, 1$               &  3  &$\scriptstyle[1:n_0+2]$    &$\scriptstyle x^{n_0}(x^{2+n_0} +6\lambda y)(x^{2+n_0}-3\lambda y)$  &$\scriptstyle\frac{n_0}{2}\frac{dx}{x}+\frac{d\ln (x^{2+n_0}+6\lambda y)}{2}+{\scriptstyle d\ln (x^{2+n_0}-3\lambda y)}$    \\
     \hline

10&$\scriptstyle1+\frac{n_0}{2}, 1$               &  3  &$\scriptstyle[1:n_0+2]$    &$\scriptstyle x^{n_0}(x^{2+n_0} +\lambda y)(x^{4+2n_0}+2\lambda x^{2+n_0}y+4\lambda^2 y^2)$  &$\frac{n_0}{2}\frac{dx}{x}+\frac{d\ln (x^{4+2n_0}+2\lambda x^{2+n_0}y+4\lambda^2 y^2)}{2}$     \\
     \hline
11&$\scriptstyle1+\frac{n_0}{2}, 1$               &  3  &$\scriptstyle[1:n_0+2]$    &$\scriptstyle x^{n_0}(x^{2+n_0} -\lambda y)(x^{2+n_0}+\lambda y)$  &$\frac{n_0}{2}\frac{dx}{x}+\frac{d\ln (x^{2+n_0}+\lambda y)}{3}$    \\
     \hline
12&$\scriptstyle1+\frac{n_0}{2}, 1$               &  3  &$\scriptstyle[1:n_0+2]$    &$\scriptstyle x^{n_0}(x^{2+n_0} +4\lambda y)(x^{2+n_0}+2\lambda y)$  &$\scriptstyle\frac{n_0}{2}\frac{dx}{x}+\frac{d\ln (x^{2+n_0}+4\lambda y)}{2}+\frac{d\ln (x^{2+n_0}+2\lambda y)}{3}$    \\
     \hline
13&$\scriptstyle1+\frac{n_0}{2}, 1$               &  3  &$\scriptstyle[1:n_0+2]$    &$\scriptstyle x^{n_0}(x^{2+n_0} -4\lambda y)(x^{4+2n_0}-2\lambda x^{2+n_0}y-2\lambda^2 y^2)$  &$\scriptstyle\frac{n_0}{2}\frac{dx}{x}+\frac{d\ln (x^{2+n_0}-4\lambda y)}{2}$     \\
     \hline
14&$\scriptstyle1+\frac{n_0}{3}, 1$               &  3  &$\scriptstyle[1:n_0+3]$    &$\scriptstyle yx^{n_0}(x^{3+n_0} -27\lambda y)$  &$\scriptstyle\frac{n_0}{3}\frac{dx}{x}+\frac{2}{3}\frac{dy}{y}+d\ln (x^{3+n_0}-27\lambda y)$     \\
     \hline

15&$\scriptstyle1+\frac{n_0}{3}, 1$               &  3  &$\scriptstyle[1:n_0+3]$    &$\scriptstyle yx^{n_0}(2x^{3+n_0} +27\lambda y)(x^{3+n_0} +54\lambda y)$  &$\scriptstyle\frac{n_0}{3}\frac{dx}{x}+\frac{2}{3}\frac{dy}{y}+\frac{d\ln (x^{3+n_0}+54\lambda y)}{2}$       \\
     \hline
16&$\scriptstyle1+\frac{n_0}{3}, 1$               &  3  &$\scriptstyle[1:n_0+3]$    &$\scriptstyle yx^{n_0}(8x^{3+n_0} +27\lambda y)$  &$\scriptstyle\frac{n_0}{3}\frac{dx}{x}+\frac{2}{3}\frac{dy}{y}$       \\
     \hline
17&$\scriptstyle1+\frac{n_0}{3}, 1$               &  3  &$\scriptstyle[1:n_0+3]$    &$\scriptstyle x^{n_0}(2x^{3+n_0}- 9\lambda y)(2x^{3+n_0}-27 \lambda y)$  &$\scriptstyle\frac{n_0}{3}\frac{dx}{x}+\frac{d\ln (2x^{3+n_0}-9\lambda y)}{2}+\frac{d\ln (2x^{3+n_0}-27\lambda y)}{2}$       \\
     \hline
18&$\scriptstyle1+\frac{n_0}{3}, 1$               &  3  &$\scriptstyle[1:n_0+3]$    &$\scriptstyle x^{n_0}(25x^{3+n_0}- 18\lambda y)(25x^{3+n_0}+27 \lambda y)$  &$\scriptstyle\frac{n_0}{3}\frac{dx}{x}+\frac{d\ln (25x^{3+n_0}-18\lambda y)}{2}$       \\
     \hline

19&$\scriptstyle1+\frac{n_0}{3}, 1$               &  3  &$\scriptstyle[1:n_0+3]$    &$\scriptstyle x^{n_0}(x^{3+n_0}- 3\lambda y)(2x^{3+n_0}+3 \lambda y)$  &$\scriptstyle\frac{n_0}{3}\frac{dx}{x}+\frac{d\ln (x^{3+n_0}-3\lambda y)}{3}$        \\
     \hline
20&$\scriptstyle1+\frac{n_0}{3}, 1$               &  3  &$\scriptstyle[1:n_0+3]$    &$\scriptstyle x^{n_0}(4x^{3+n_0}- 75\lambda y)(4x^{3+n_0}+15\lambda y)(x^{3+n_0}-30\lambda y)$  &$\scriptstyle\frac{n_0}{3}\frac{dx}{x}+\frac{d\ln (x^{3+n_0}-30\lambda y)}{2}+\frac{d\ln (4x^{3+n_0}-75\lambda y)}{3}$        \\
     \hline

21&$\scriptstyle1+\frac{n_0}{3}, 1$               &  3  &$\scriptstyle[1:n_0+3]$    &$\scriptstyle x^{n_0}(x^{3+n_0}+27\lambda y)(x^{3+n_0}-9\lambda y)$  &$\scriptstyle\frac{n_0}{3}\frac{dx}{x}+\frac{d\ln (x^{3+n_0}+27\lambda y)}{2}+\frac{3}{2}d\ln (x^{3+n_0}-9\lambda y)$        \\
     \hline

22&$\scriptstyle1+\frac{n_0}{3}, 1$               &  3  &$\scriptstyle[1:n_0+3]$    &$\scriptstyle x^{n_0}(x^{3+n_0}-18\lambda y)(2x^{3+n_0}-27\lambda y)(4x^{3+n_0}+9\lambda y)$  &$\scriptstyle \frac{n_0}{3}\frac{dx}{x}+\frac{d\ln (x^{3+n_0}-18\lambda y)}{2}$        \\
     \hline
23&$\scriptstyle1+\frac{n_0}{3}, 1$               &  3  &$\scriptstyle[1:n_0+3]$    &$\scriptstyle x^{n_0}(8x^{3+n_0}-27\lambda y)(8x^{3+n_0}+9\lambda y)$  &$\scriptstyle\frac{n_0}{3}\frac{dx}{x}$        \\
     \hline

 \end{tabular}

\hspace{-1.0cm}\begin{tabular}{|c|c|c|c|c|c|}
     \hline

24&$\scriptstyle1+\frac{n_0}{3}, 1$               &  3  &$\scriptstyle[1:n_0+3]$    &$\scriptstyle x^{n_0}(x^{3+n_0}+9\lambda y)(x^{6+2n_0}+36\lambda x^{3+n_0}y+972\lambda^2 y^2)$  &$\scriptstyle\frac{n_0}{3}\frac{dx}{x}+\frac{d\ln (x^{6+2n_0}+36\lambda x^{3+n_0}y+972\lambda^2 y^2)}{2}$        \\
     \hline
25&$\scriptstyle1+\frac{n_0}{3}, 1$               &  3  &$\scriptstyle[1:n_0+3]$    &$\scriptstyle x^{n_0}(x^{3+n_0}-216\lambda y)(x^{3+n_0}-144\lambda y)$  &$\scriptstyle\frac{n_0}{3}\frac{dx}{x}+\frac{3}{2}d\ln (x^{3+n_0}-144\lambda y)$        \\
     \hline
26&$\scriptstyle1+\frac{n_0}{3}, 1$               &  3  &$\scriptstyle[1:n_0+3]$    &$\scriptstyle x^{n_0}(25x^{3+n_0}+432\lambda y)(25x^{3+n_0}+72\lambda y)$  &$\scriptstyle \frac{n_0}{3}\frac{dx}{x}+\frac{d\ln (25x^{3+n_0}+432\lambda y)}{2}$        \\
     \hline
\end{tabular}\\

\bigskip

\end{document}